\documentclass[11pt,reqno]{amsart}

\usepackage{amsmath,amssymb,mathrsfs,amsthm}
\usepackage{graphicx,cite,cases}
\usepackage{xcolor}
\newtheorem{theorem}{Theorem}[section]

\newtheorem{lemma}[theorem]{Lemma}

\numberwithin{equation}{section}

\begin{document}

\title[the PML problem for the biharmonic wave scattering]{Numerical solution to the PML problem of the biharmonic wave scattering in periodic structures}

\author{Peijun Li}
\address{LSEC, ICMSEC, Academy of Mathematics and Systems Science, Chinese Academy of Sciences, Beijing 100190, China, and School of Mathematical Sciences, University of Chinese Academy of Sciences, Beijing 100049, China}
\email{lipeijun@lsec.cc.ac.cn}

\author{Xiaokai Yuan}
\address{School of Mathematics, Jilin University, Changchun 130012, Jilin, China}
\email{yuanxk@jlu.edu.cn}

\thanks{The second
author is supported by the NSFC grants 12201245 and 12171017.}

\subjclass[2010]{78A45, 65N30}

\keywords{Biharmonic wave equation, cavity scattering problem, PML method, convergence analysis}

\maketitle

\begin{abstract}
Consider the interaction of biharmonic waves with a periodic array of cavities, characterized by the Kirchhoff--Love model. This paper investigates the perfectly matched layer (PML) formulation and its numerical soution to the governing biharmonic wave equation. The study establishes the well-posedness of the associated variational problem employing the Fredholm alternative theorem. Based on the examination of an auxiliary problem in the PML layer, exponential convergence of the PML solution is attained. Moreover, it develops and compares three decomposition methods alongside their corresponding mixed finite element formulations, incorporating interior penalty techniques for solving the PML problem. Numerical experiments validate the effectiveness of the proposed methods in absorbing outgoing waves within the PML layers and suppressing oscillations in the bending moment of biharmonic waves near the cavity's surface.
\end{abstract}

\section{introduction}\label{Section1}

The biharmonic wave equation governs the propagation of flexural waves and plays an important role in thin plate elasticity, including applications such as ultra-broadband elastic cloaking \cite{DZAL-PRL-2018, FGE-PRL-2009, FGE-JCP-2011}, platonic crystals \cite{S-UA-2013, H-UL-2014, HCMM-RRSA-2016}, and acoustic black holes \cite{PGCS-JSV-2020}. However, most literature primarily focuses on addressing the bi-Laplacian or biharmonic problem, which concerns the static state of the biharmonic wave equation, in a bounded domain. Because of the intricate nature of constructing $C^1$ continuous piecewise polynomials on meshes, employing standard $H^2$ conforming finite elements for fourth-order equations becomes less appealing. Consequently, various nonconforming and discontinuous finite element methods have been proposed, including the weak Galerkin finite element method  \cite{MWZ-JSC-2014, YZ-SIAM-2020, ZZ-JSC-2015}, the virtual element method \cite{AMV-MMMAS-2018, ZCZ-MMMAS-2016}, and the mixed element method \cite{AD-NM-2001, C-FEM-1978, CR-1974, CG-CMAME-1975}.

Fewer mathematical and computational results have been reported for the biharmonic wave equation, especially in unbounded domains, due to additional complexities inherent in wave phenomena. In \cite{DL-2023-arXiv}, a novel boundary integral equation formulation was developed for the biharmonic wave scattering problem involving bounded cavities. By introducing two auxiliary functions, the biharmonic wave equation is decomposed into the coupled Helmholtz and modified Helmholtz equations. The boundary integral equation method is employed to transform equivalently the scattering problem from an unbounded domain onto the surface of the cavity. The well-posedness of the coupled system is established using the Riesz--Fredholm theorem. In \cite{BCF-CMS-2019, BH-SIAM-2020}, the variational approach was empolyed, and the transparent boundary condition (TBC) was deduced on an artificial boundary enclosing the cavity, thereby converting the unbounded domain into a bounded one. The well-posedness was demonstrated for the corresponding variational formulation of the cavity scattering problem with various boundary conditions. Numerically, by decomposing the biharmonic wave equation into the coupled Helmholtz and modified Helmholtz equations, the equivalent TBC was constructed in \cite{YLYZ-RAM-2023, YL-JCP-2023} for each wave component, and a mixed finite element method was proposed by incorporating interior and boundary penalty schemes to suppress the oscillations of the bending moment near the cavity's surface.

Ever since  B\'erenger introduced the perfectly matched layer (PML) method in \cite{B-JCP-1994}, it has emerged as a widely adopted domain truncation technique for simulating diverse wave scattering problems in unbounded domains. The fundamental concept is to enclose the domain of interest with a sufficiently large artificial absorbing layer, followed by truncating the computational domain with either Dirichlet or Neumann boundary condition applied on the outer boundary of the PML layer. In comparison to the TBC method, which requires handling nonlocal Dirichlet-to-Neumann (DtN) operators, the PML method is more convenient to implement due to its easier management of local boundary conditions. In \cite{BLY-2023-sub}, the PML method was investigated for the biharmonic wave scattering problem in periodic structures. 
The equivalent PML-DtN operator was derived, and then the sesquilinear form was verified to satisfy G\r{a}rding's inequality. Subsequently, the well-posedness of the variational problem was established by applying the Fredholm alternative theorem. The error analysis between the solution of the PML problem and that of the original scattering problem was conducted by directly estimating the difference between the PML-DtN and the original DtN operators. This approach has been adopted to analyze the error of the PML method for various wave scattering problems, including acoustic and electromagnetic wave scattering problems \cite{CW-SINUM-2003, LWZ-SIMA-2011, BW-SINUM }. However, in cases where the scattering problem is complex and the expression of the TBC is intricate, validating whether the PML-DtN operator is well-defined can become more challenging, if not impossible.
 
In this paper, we investigate the PML formulation and its numerical solution for the scattering of biharmonic waves by 
cavities embedded in an infinitely extending elastic thin plate. The work comprises three main contributions:

\begin{enumerate}
 
 \item[(1)] Establishment of the well-posedness for the PML problem.
 
 \item[(2)] Achievement of exponential convergence for the PML solution.
 
 \item[(3)] Proposal of mixed finite element methods for the PML problem. 
 
\end{enumerate}

Specifically, we investigate the scattering behavior of a plane incident wave as it interacts with a periodic array of one-dimensional cavities, described by the biharmonic wave equation imposed in an unbounded domain. Through the application of complex coordinate stretching, we introduce the PML formulation to convert the scattering problem from an unbounded domain to a bounded domain. In \cite{BLY-2023-sub}, the well-posedness was established for the associated variational problem by using the PML-DtN operator. In this work, we present an alternative approach to establish the well-posedness of the PML problem by directly studying the associated variational problem in the computational domain. Regarding the error analysis, we adopt the methodology outlined in \cite{CL-SINUM-2005, BP-MC-2007, BPT-MC-2010}, where acoustic, electromagnetic, and elastic wave scattering problems were addressed. The objective remains to analyze the difference between the PML-DtN and original DtN operators. Instead of explicitly generating the PML-DtN operator and computing the difference, an auxiliary problem within the PML layer is introduced. In this auxiliary problem, one boundary condition is defined by the discrepancy between the two operators, while the other boundary condition is defined by the propagating operator.
This approach is more convenient for demonstrating the exponential convergence of the propagating operator with respect to the thickness of the layer. The difference between the DtN operators can then be bounded by the propagating operator using the stability estimate of the auxiliary problem. Finally, the error estimate of the PML solution can be obtained.
 
Furthermore, we propose three mixed finite element methods incorporating the interior penalty technique. These methods are based on three distinct decompositions of the biharmonic wave equation. Numerical experiments demonstrate that the proposed methods effectively absorb outgoing waves within the PML layers and suppress oscillations of the bending moment near the surface of the cavities.

The structure of the paper is outlined as follows. Section \ref{Section2} introduces the model equation. The well-posedness of the scattering problem is discussed in Section \ref{Section3}. Section \ref{section4} employs the PML method to reformulate the problem from an unbounded domain to a bounded computational domain, and addresses the well-posedness of the associated variational problem. Convergence analysis is detailed in Section \ref{section5}. Section \ref{section6} presents the mixed finite element methods along with accompanying numerical experiments. Finally, general remarks are provided in Section \ref{section7} as the paper concludes.

\section{Problem formulation}\label{Section2}

This paper investigates the scattering of a plane incident wave by cavities, which is characterized by the Kirchhoff--Love model. The cavities are periodically arranged along the $x_1$-axis within an infinitely extending elastic thin plate, with the periodicity of alignment assumed to be $\Lambda$. 

\begin{figure}
\centering
\includegraphics[width=0.6\textwidth]{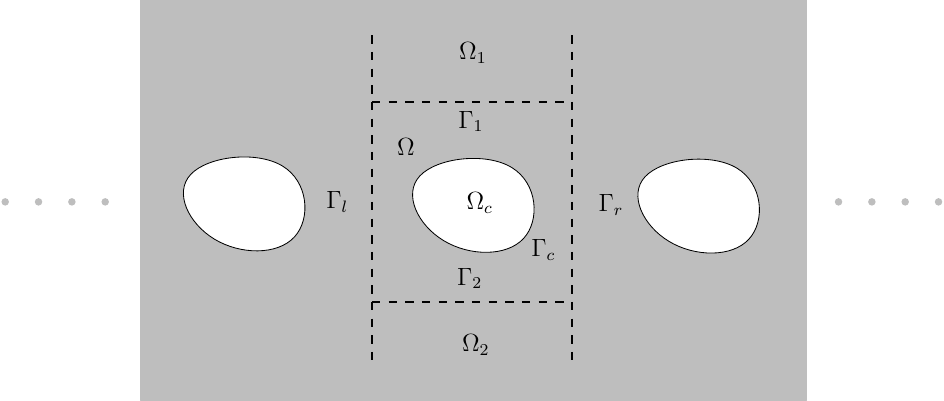}
\caption{Schematic of the problem geometry.}
\label{pg}
\end{figure}

Given the periodic nature of the structure, it enables us to confine the problem to just one periodic unit. Denote by $\Omega_c$ the region of the cavity, which is assumed to have a Lipschitz continuous boundary $\Gamma_c$. Let $R$ be a rectangle of large size enclosing the cavity $\Omega_c$. For the sake of simplicity, we define $R$ as follows: $R=\left\{x\in\mathbb{R}^2: 0<x_1<\Lambda, h_2<x_2<h_1\right\}$, where $h_k, k=1, 2$ are constants. Let $\Gamma_k$ be defined as $\left\{x\in\mathbb{R}^2: 0<x_1<\Lambda, x_2=h_k\right\}$ for $k=1, 2$, $\Gamma_l$ as $\left\{x\in\mathbb{R}^2: x_1=0, h_2<x_2<h_1\right\}$, and $\Gamma_r$ as $\left\{x\in\mathbb{R}^2: x_1=\Lambda, h_2<x_2<h_1\right\}$. Consider $\Omega$ as the set difference $R\setminus\overline{\Omega_c}$. Define $\Omega_1$ and $\Omega_2$ as the areas situated above $\Gamma_1$ and below $\Gamma_2$, respectively. The problem geometry is illustrated in Figure \ref{pg}. 

The total field $u$ satisfies the biharmonic wave equation
\begin{equation}\label{TotalBiharmonic}
 \Delta^2 u-\kappa^4 u=0\quad \text{in} ~ \Omega, 
\end{equation}
where $\kappa>0$ is the wavenumber, and $u$ denotes the out-of-plane displacement of the plate. We consider the clamped boundary condition on $\Gamma_c$ without any loss of generality:
\begin{equation}\label{cbc}
 u=0,\quad \partial_{\nu} u=0, 
\end{equation}
where $\nu$ stands for the unit normal vector on $\Gamma_c$.

Consider a time-harmonic plane wave as the incident wave, represented by
\[
u^i(x)=e^{{\rm i}(\alpha x_1-\beta x_2)}, \quad x\in\mathbb R^2,
\]
where $\alpha = \kappa\sin\theta$ and $\beta = \kappa\cos\theta$, and $\theta$ belongs to the interval $\left(-\frac{\pi}{2}, \frac{\pi}{2}\right)$, representing the incident angle. It is clear to note that the incident field $u^i$ satisfies 
\[
 \Delta^2 u^i - \kappa^4 u^i=0\quad\text{in} ~ \mathbb R^2.
\]

Due to the periodic nature of the structure and the incident wave, the solution to \eqref{TotalBiharmonic}--\eqref{cbc} is quasi-periodic. Precisely, if $u$ is designated as a quasi-periodic function,  then $u(x) e^{-{\rm i}\alpha x_1}$ becomes a periodic function of $x_1$ with a period of $\Lambda$. The quasi-periodic nature of $u$ enforces a quasi-periodic boundary condition on $\Gamma_l$ and $\Gamma_r$, expressed as $u(0, x_2) = e^{-{\rm i}\alpha\Lambda}u(\Lambda, x_2)$. Additionally, the radiation condition needs to be imposed: the scattered field $u^s=u-u^i$ within $\Omega_1$ and the total field $u$ within $\Omega_2$ are bounded outgoing waves.

Denote by $H^2(\Omega)$ the standard Sobolev space. Introduce the function space with quasi-periodicity:
\[
H_{\rm qp}^2(\Omega)=\left\{u\in H^2(\Omega): u(\Lambda, x_2)=u(0, x_2)e^{{\rm i}\alpha\Lambda}\right\},
\]
and its subspace
\[
H_{{\rm qp}, \Gamma_c}^{2}(\Omega)=\left\{u\in H^2_{\rm qp}(\Omega): u=0, \partial_{\nu} u=0\text{ on }\Gamma_c\right\}.
\]
It is evident that $H_{\rm qp}^2(\Omega)$ and $H_{{\rm qp}, \Gamma_c}^{2}(\Omega)$ are subspaces of $H^2(\Omega)$, utilizing the conventional $H^2$-norm.

For any $u\in H_{\rm qp}^2(\Omega)$, it admits a Fourier series expansion over $\Gamma_k, k=1, 2:$
\[
u(x, h_k)=\sum\limits_{n\in\mathbb{Z}} u^{(n)}(h_k) e^{{\rm i}\alpha_n x_1},
\]
where
\[
 \alpha_n=\alpha+n\left(\frac{2\pi}{\Lambda}\right), \quad u^{(n)}(h_k)=\frac{1}{\Lambda}\int_{0}^{\Lambda} 
 u(x, h_k) e^{-{\rm i}\alpha_n x_1}{\rm d}x_1. 
\]
The trace function space $H^s(\Gamma_k)$, where $s \in \mathbb{R}$, is defined as follows:
\[
H^{s}(\Gamma_k)=\left\{u\in L^2(\Gamma_k): \|u\|_{H^s(\Gamma_k)}<\infty\right\},
\]
with the norm given by
\[
\|u\|_{H^s(\Gamma_k)}=\left(\Lambda\sum\limits_{n\in\mathbb{Z}}\left(1+\alpha_n^2\right)^s |u^{(n)}(h_k)|^2\right)^{1/2}.
\]

In this paper, the notation $a\lesssim b$ denotes $a\leq Cb$ with $C>0$ being a constant. The specific value of $C$ is not essential, but its dependence is evident within the given context.

\section{The reduced problem}\label{Section3}

The original scattering problem \eqref{TotalBiharmonic}--\eqref{cbc} is formulated in an unbounded domain. To facilitate computation, it is preferable to equivalently transform the problem into a bounded domain. In \cite{BLY-2023-sub}, it is demonstrated that \eqref{TotalBiharmonic}--\eqref{cbc} is equivalent to the following boundary value problem, incorporating transparent boundary conditions on $\Gamma_1\cup\Gamma_2:$
\begin{equation}\label{TotalBiharmonicReduced}
\left\{
\begin{aligned}
	& \Delta^2 u-\kappa^4 u=0 &&\text{in} ~ \Omega,\\
	& u=0,\quad \partial_{\nu} u=0 & &\text{on} ~ \Gamma_c,\\
	& N_1 u=T_{11}^{(1)} f_1+T_{12}^{(1)}g_1+p_1 &&\text{on} ~ \Gamma_1,\\
	& M_1 u=T_{21}^{(1)}f_1+T_{22}^{(1)}g_1+p_2 &&\text{on} ~ \Gamma_1,\\
	& N_2 u=T_{11}^{(2)} f_2+T_{12}^{(2)}g_2 & &\text{on} ~ \Gamma_2,\\
	& M_2 u=T_{21}^{(2)}f_2+T_{22}^{(2)}g_2 & &\text{on} ~ \Gamma_2,
\end{aligned}
\right.
\end{equation}
where the boundary differential operators $M_k, N_k$  on $\Gamma_k,  k=1, 2$ are given by  (cf. \cite{HW-2021})
\begin{eqnarray}
	&&N_1 u=-(2-\mu)\frac{\partial^3 u}{\partial x_1^2\partial x_2}-\frac{\partial^3 u}{\partial x_2^3},\quad
	M_1 u=\mu\frac{\partial^2 u}{\partial x_1^2}+\frac{\partial^2 u}{\partial x_2^2}, \label{MNupper}\\
	&&N_2 u=(2-\mu)\frac{\partial^3 u}{\partial x_1^2\partial x_2}+\frac{\partial^3 u}{\partial x_2^3},\quad
	M_2 u=\mu\frac{\partial^2 u}{\partial x_1^2}+\frac{\partial^2 u}{\partial x_2^2}, \label{MNlower}	
\end{eqnarray}
where $\mu\in(0, 1)$ is the Poisson ratio, and the DtN operators $T^{(k)}_{i,j}, k=1, 2, i,j=1,2$ are given by 
\begin{equation}\label{T1}
\left\{
\begin{aligned}
	(T_{11}^{(k)} f)(x_1) &=\sum\limits_{n\in\mathbb{Z}}
	{\rm i}\beta_n\gamma_n\left(\gamma_n-{\rm i}\beta_n\right)f^{(n)} e^{{\rm i}\alpha_n x_1},\\
	(T_{21}^{(k)} f)(x_1) &=-\sum\limits_{n\in\mathbb{Z}} \left(\mu\alpha_n^2
	-{\rm i}\beta_n\gamma_n\right) f^{(n)} e^{{\rm i}\alpha_n x_1},\\
	(T_{12}^{(k)}g)(x_1) &=-\sum\limits_{n\in\mathbb{Z}} \left(\mu\alpha_n^2-{\rm i}\beta_n\gamma_n\right) g^{(n)} 
	e^{{\rm i}\alpha_n x_1},\\
	(T_{22}^{(k)}g)(x_1) &=-\sum\limits_{n\in\mathbb{Z}} \left(\gamma_n-{\rm i}\beta_n\right)g^{(n)} e^{{\rm i}\alpha_n x_1},
\end{aligned}
\right.
\end{equation}
where $(f_k, g_k)=(u, \partial_{\nu} u)|_{\Gamma_k}, k=1, 2$, stand for the Dirichlet data on $\Gamma_k$, respectively.  
The inhomogenous functions $p_1, p_2$ are given by
\begin{eqnarray}\label{alphabeta}
&&p_1(x_1)=-\left(2{\rm i}\beta\alpha^2+2\beta^2\gamma\right) e^{{\rm i}(\alpha x_1-\beta h_1)},\quad
p_2(x_1)=-(2\beta^2+2{\rm i}\beta\gamma)e^{{\rm i}(\alpha x_1-\beta h_1)},\label{p12}
\end{eqnarray}
where
\begin{eqnarray*}\label{betagamman}
	\beta_n=\left\{
	\begin{aligned}
	&(\kappa^2-\alpha_n^2)^{1/2}\quad & &\text{if} ~ \kappa>|\alpha_n|,\\
	&{\rm i}(\alpha_n^2-\kappa^2)^{1/2}\quad& &\text{if} ~ \kappa<|\alpha_n|,
	\end{aligned}
	\right. \quad
	\gamma_n=(\kappa^2+\alpha_n^2)^{1/2},\quad
	\gamma=\gamma_0.
\end{eqnarray*}

The variational problem of \eqref{TotalBiharmonicReduced} is to find 
 $u\in H_{{\rm qp}, \Gamma_c}^2(\Omega)$ such that 
\begin{equation}\label{variationalTBC}
	a(u, v)=\int_{\Gamma_1} \left(p_1\bar{v} + p_2 \partial_\nu \bar{v}\right){\rm d}s,
	\quad \forall\, v\in H_{{\rm qp}, \Gamma_c}^2(\Omega),
\end{equation}
where the sesquilinear form $a(u, v):  H_{{\rm qp}, \Gamma_c}^2(\Omega)\times  
 H_{{\rm qp}, \Gamma_c}^2(\Omega)\rightarrow \mathbb{C}$ is defined as 
\begin{align}\label{tbcauv}
	a(u, v)=\int_{\Omega} \bigg[\mu\Delta u\Delta \bar{v}+(1-\mu)\sum\limits_{i,j=1}^2
	 \frac{\partial^2 u}{\partial x_i\partial x_j}\frac{\partial^2 \bar{v}}{\partial x_i\partial x_j}
	 -\kappa^4 u\bar{v}\bigg]{\rm d}x -\sum_{k=1}^2 \int_{\Gamma_k}  (\mathbb T^{(k)}
	 \boldsymbol{u})\cdot\overline{\boldsymbol v}{\rm d}s,
\end{align}
with $\boldsymbol{u}, \boldsymbol{v}$, and $\mathbb T^{(k)}$ given by 
\[
 \boldsymbol{u}=\begin{bmatrix}
                 u\\
                 \partial_\nu u
                \end{bmatrix},\quad \boldsymbol{v}=\begin{bmatrix}
                 v\\
                 \partial_\nu v
                \end{bmatrix},\quad 
\mathbb T^{(k)}=\begin{bmatrix}
                 T_{11}^{(k)} & T_{12}^{(k)}\\[5pt]
                 T_{21}^{(k)} & T_{22}^{(k)}
                \end{bmatrix}.
\]

The properties of the DtN operators \eqref{T1}, as well as the well-posedness of the variational problem \eqref{variationalTBC}, are discussed in \cite{BLY-2023-sub}. The well-posedness is presented below.

\begin{theorem}\label{MainResult1}
There exists a unique weak solution $u\in H^2_{{\rm qp}, \Gamma_c}(\Omega)$ to the variational problem \eqref{variationalTBC}, with the exception of a discrete set of wavenumbers $\kappa$.
\end{theorem}

\section{The PML problem}\label{section4}

In this section, we introduce the PML formulation of the scattering problem \eqref{TotalBiharmonic}--\eqref{cbc} and examine the well-posedness of its variational problem.

\begin{figure}
\centering
\includegraphics[width=0.22\textwidth]{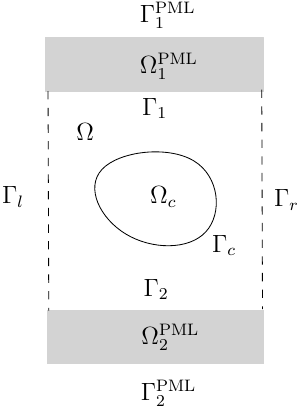}
\caption{The schematic of the PML problem.}
\label{pml}
\end{figure}

\subsection{ The PML formulation}

Let $\Omega_k^{\rm PML}$, where $k=1, 2$, represent the PML regions located above and below the interfaces $\Gamma_k$, each with a thickness denoted by $\Delta h_k$. Let $\Gamma_k^{\rm PML}$ denote the outer boundaries of $\Omega_k^{\rm PML}$, $k=1, 2$. Let $\Omega^{\rm PML}:=\Omega\cup\Omega_1^{\rm PML}\cup\Omega_2^{\rm PML}$ be the computational domain where the PML problem is defined. Figure \ref{pml} illustrates a schematic of the PML problem.

The PML formulation is introduced through complex coordinate stretching (cf. \cite{TC-IEEE-1997}):
\begin{equation}\label{pmlcons}
\tilde{x}_2=\varphi(x_2)=\int_0^{x_2} \sigma(t){\rm d}t,
\end{equation}
where 
\[
\sigma(t)=\left\{
\begin{aligned}
& 1 & &\text{if} ~ h_2\leq t\leq h_1,\\
& 1+\sigma_1\left(\frac{t-h_1}{\Delta h_1}\right)^m+{\rm i}\sigma_2\left(\frac{t-h_1}{\Delta h_1}\right)^m & &\text{if} ~  t>h_1,\\
& 1+\sigma_1\left(\frac{h_2-t}{\Delta h_2}\right)^m+{\rm i}\sigma_2\left(\frac{h_2-t}{\Delta h_2}\right)^m & &\text{if} ~ t<h_2,
\end{aligned}
\right.
\]
with $m$ being required to be greater than $3$ and the PML parameters $\sigma_1$ and  $\sigma_2$ being positive constants satifying the condition
\begin{equation}\label{PMLrestriction}
	1+\sigma_1>\sqrt{\frac{3+\mu}{1-\mu}}\sigma_2,
\end{equation}
which is required in Subsections \ref{PML:variational} and \ref{PML:wellp}.

Let $\tilde{u}$ be the solution of the total field for the truncated PML problem, which satisfies
\begin{equation}\label{MBHE}
\left\{
\begin{aligned}
&\tilde{\Delta}^2\tilde{u}-\kappa^4\tilde{u}=f & &{\rm in} ~ \Omega^{\rm PML},\\
&\tilde{u}=u^i,\quad\partial_{\nu}\tilde{u}=\partial_{\nu} u^i & &{\rm on} ~ \Gamma_1^{\rm PML},\\
&\tilde{u}=0,\quad\partial_{\nu}\tilde{u}=0 & &{\rm on} ~ \Gamma_2^{\rm PML},\\
&\tilde{u}=0,\quad\partial_{\nu}\tilde{u}=0&&{\rm on}~\Gamma_c,
\end{aligned}
\right.
\end{equation}
where $\tilde x=(x_1, \tilde x_2)$ represents the complex coordinates and $\tilde{\Delta}$ is the Laplacian operator in the complex coordinates, given by 
\[
\tilde{\Delta}\tilde{u}=\frac{\partial^2 \tilde{u}}{\partial x_1^2}+\frac{\partial^2 \tilde{u}}{\partial \tilde{x}_2^2},
\]
and 
\begin{equation}\label{fpml}
f(\tilde x)=\left\{
\begin{aligned}
& \tilde{\Delta}^2 u^i-\kappa^4 u^i & &\text{in} ~ \Omega_1^{\rm PML},\\
& 0  & &\text{in} ~ \Omega\setminus\overline{\Omega_1^{\rm PML}}. 
\end{aligned}
\right.
\end{equation}

In this context, we regard the complex coordinate variable $\tilde x=(x_1, \tilde{x}_2)$ as being in $\Omega^{\rm PML}_i$ or situated on $\Gamma_i^{\rm PML}$ if $x=(x_1, \psi(\tilde{x}_2))$ is in $\Omega^{\rm PML}_i$ or on $\Gamma_i^{\rm PML}$,  $i=1, 2$, respectively, and $\psi$ is the inverse function of $\varphi$, i.e., $x_2=\psi(\tilde x_2)=\varphi^{-1}(\tilde x_2)$.

\subsection{The variational problem}\label{PML:variational}

Denote by $\hat{u}(x_1, \psi(\tilde{x}_2))=\tilde{u}(x_1, \tilde{x}_2)$ and $\hat{f}(x_1, \psi(\tilde{x}_2))=f(x_1, \tilde{x}_2)$. It follows from the change of variables that
\begin{align}\label{MBHE2}
	\tilde{\Delta}^2 \tilde{u}-\kappa^4\tilde{u}
	&=\frac{\partial^2}{\partial x_1^2}\left(\frac{\partial^2 \hat{u}}{\partial x_1^2}
	+\mu\frac{1}{\sigma}\frac{\partial}{\partial x_2}\left(\frac{1}{\sigma}\frac{\partial \hat{u}}{\partial x_2}\right)\right)+2(1-\mu)\frac{1}{\sigma}\frac{\partial^2}{\partial x_1 \partial x_2}
		\left(\frac{1}{\sigma}\frac{\partial^2 \hat{u}}{\partial x_1 \partial x_2}\right)\notag\\
	&\quad +\frac{1}{\sigma}\frac{\partial}{\partial x_2}\left(
	\frac{1}{\sigma}\frac{\partial}{\partial x_2}\left(\frac{1}{\sigma}\frac{\partial}{\partial x_2}
		\left(\frac{1}{\sigma}\frac{\partial \hat{u}}{\partial x_2}\right)
		+\mu\frac{\partial^2 \hat{u}}{\partial x_1^2}\right)\right)-\kappa^4 \hat{u}.
\end{align}
For convenience, we do not distinguish between $\hat{u}$, $\hat{f}$, and $\tilde{u}$, $f$ in the rest of the paper when there is no confusion.

Introduce two subspaces of $H^2_{{\rm qp}, \Gamma_c}(\Omega^{\rm PML})$:
\begin{align*}
H^2_{{\rm qp}, 0}(\Omega^{\rm PML})& :=\left\{u\in H_{{\rm qp}, \Gamma_c}^2(\Omega^{\rm PML}):
u=0,\,\partial_\nu u=0\,{\rm on}\,\,\Gamma^{\rm PML}_2\right\},\\
H^2_{{\rm qp}, 00}(\Omega^{\rm PML})& :=\left\{u\in H_{{\rm qp}, \Gamma_c}^2(\Omega^{\rm PML}): 
u=0,\,\partial_\nu u=0\,{\rm on}\,\,\Gamma^{\rm PML}_2\cup\Gamma^{\rm PML}_{1}\right\}.
\end{align*}
Multiplying both sides of \eqref{MBHE2} by $v\in H_{{\rm qp}, 00}^2(\Omega^{\rm PML})$ yields 
\begin{align}
	\int_{\Omega^{\rm PML}} \sigma f\overline{v}\,{\rm d}x 
	&= \int_{\Omega\cup\Omega_1^{\rm PML}\cup\Omega_2^{\rm PML}} 
	\Bigg[\frac{\partial^2}{\partial x^2_1}\left(\sigma\frac{\partial^2 \tilde{u}}{\partial x_1^2}
	+\mu\frac{\partial}{\partial x_2}\left(\frac{1}{\sigma}\frac{\partial \tilde{u}}{\partial x_2}\right)\right)\notag\\
	&\quad +\frac{\partial}{\partial x_2}\left(
	\frac{1}{\sigma}\frac{\partial}{\partial x_2}\left(\frac{1}{\sigma}\frac{\partial}{\partial x_2}
		\left(\frac{1}{\sigma}\frac{\partial \tilde{u}}{\partial x_2}\right)+\mu\frac{\partial^2 \tilde{u}}{\partial x_1^2}\right)\right)\notag\\
    &\quad +2(1-\mu)\frac{\partial^2}{\partial x_1 \partial x_2}\left(\frac{1}{\sigma}\frac{\partial^2 \tilde{u}}{\partial x_1 \partial x_2}\right)-\kappa^4 \sigma \tilde{u}\Bigg]\overline{v}{\rm d}x\label{auv}.
\end{align}
Since $\sigma=1$ in $\Omega$, we obtain from the Green's formula that
\begin{align}\label{Physical}
 \int_{\Omega} \Bigg[\frac{\partial^2}{\partial x^2_1}\left(\frac{\partial^2 \tilde{u}}{\partial x_1^2}
	+\mu\frac{\partial^2 \tilde{u}}{\partial x_2^2}\right)
	+2(1-\mu)\frac{\partial^2}{\partial x_1 \partial x_2}\left(\frac{\partial^2 \tilde{u}}{\partial x_1 \partial x_2}\right)
	+\frac{\partial^2}{\partial x_2^2}\left(
		\frac{\partial^2 \tilde{u}}{\partial x^2_2}+\mu\frac{\partial^2 \tilde{u}}{\partial x_1^2}\right)
		-\kappa^4\tilde{u}\Bigg]\overline{v}{\rm d}x\notag\\
=\int_{\Omega} \left[\mu\Delta \tilde{u}\Delta\overline{v}
	+(1-\mu)\sum\limits_{i, j=1}^2\frac{\partial^2 \tilde{u}}{\partial x_i\partial x_j}
\frac{\partial^2 \overline{v}}{\partial x_i\partial x_j}-\kappa^4 \tilde{u}\overline{v}\right]{\rm d}x-\sum\limits_{k=1}^2\int_{\Gamma_k}\left(\overline{v} N_k \tilde{u}+\partial_{\nu} \overline{v} M_k \tilde{u}\right){\rm d}s,
\end{align}
where $M_k, N_k, k=1, 2$ are given in \eqref{MNupper}--\eqref{MNlower}. In $\Omega^{\rm PML}_1$, it follows from the quasi-periodic boundary condition and Green's formula that
\begin{align*}
& \int_{\Omega^{\rm PML}_1} \overline{v}\,\frac{\partial^2}{\partial x_1^2}\left(\sigma\frac{\partial^2 \tilde{u}}{\partial x_1^2}
	+\mu\frac{\partial}{\partial x_2}\left(\frac{1}{\sigma}\frac{\partial\tilde{u}}{\partial x_2}\right)\right)
	{\rm d}x\\
&=\int_{h_1}^{h_1+\delta h_1} \overline{v}\,\frac{\partial}{\partial x_1}\left(\sigma\frac{\partial^2 \tilde{u}}{\partial x_1^2}
	+\mu\frac{\partial}{\partial x_2}\left(\frac{1}{\sigma}\frac{\partial\tilde{u}}{\partial x_2}\right)\right)\Big|_{0}^{\Lambda}{\rm d}x_2
	-\int_{\Omega^{\rm PML}_1} \frac{\partial \overline{v}}{\partial x_1}\,\frac{\partial}{\partial x_1}\left(\sigma\frac{\partial^2 \tilde{u}}{\partial x_1^2}
	+\mu\frac{\partial}{\partial x_2}\left(\frac{1}{\sigma}\frac{\partial\tilde{u}}{\partial x_2}\right)\right)
	{\rm d}x\\
&=-\int_{h_1}^{h_1+\delta h_1}\frac{\partial \overline{v}}{\partial x_1}\,
	\left(\sigma\frac{\partial^2 \tilde{u} }{\partial x_1^2}
	+\mu\frac{\partial}{\partial x_2}\left(\frac{1}{\sigma}\frac{\partial\tilde{u}}{\partial x_2}\right)\right)\Big|_{0}^{\Lambda}{\rm d}x_2
	+\int_{\Omega^{\rm PML}_1} \frac{\partial^2 \overline{v}}{\partial x_1^2}\,\left(\sigma\frac{\partial^2 \tilde{u}}{\partial x_1^2}
	+\mu\frac{\partial}{\partial x_2}\left(\frac{1}{\sigma}\frac{\partial\tilde{u}}{\partial x_2}\right)\right)
	{\rm d}x\\
&=	\int_{\Omega^{\rm PML}_1} \frac{\partial^2 \overline{v}}{\partial x_1^2}\,\left(\sigma\frac{\partial^2 \tilde{u}}{\partial x_1^2}
	+\mu\frac{\partial}{\partial x_2}\left(\frac{1}{\sigma}\frac{\partial\tilde{u}}{\partial x_2}\right)\right)
	{\rm d}x.
\end{align*}

Similarly, by applying the Green's formula and the Dirichlet boundary condition on $\Gamma_1^{\rm PML}$, we obtain 
\begin{align*}
\int_{\Omega^{\rm PML}_1} 
	\overline{v}\,\frac{\partial^2}{\partial x_1 \partial x_2}
	\left(\frac{1}{\sigma}\frac{\partial^2 \tilde{u}}{\partial x_1 \partial x_2}\right){\rm d}x
=\int_{\Omega^{\rm PML}_1}\frac{1}{\sigma}
\frac{\partial^2\overline{v}}{\partial x_1 \partial x_2}\,\frac{\partial^2 \tilde{u}}{\partial x_1 \partial x_2}{\rm d}x
-\int_{\Gamma_1}\overline{v}\,\frac{1}{\sigma}\frac{\partial^3 \tilde{u}}{\partial x_1^2 \partial x_2}{\rm d}s
\end{align*}
and
\begin{align*}
&\int_{\Omega^{\rm PML}_1} \overline{v}\,\frac{\partial}{\partial x_2}\left(
	\frac{1}{\sigma}\frac{\partial}{\partial x_2}\left(\frac{1}{\sigma}\frac{\partial}{\partial x_2}
		\left(\frac{1}{\sigma}\frac{\partial \tilde{u}}{\partial x_2}\right)+\mu\frac{\partial^2 \tilde{u}}{\partial x_1^2}\right)\right)
		{\rm d}x\\
&=	-\int_{\Gamma_1} \overline{v}\,	\left(
	\frac{1}{\sigma}\frac{\partial}{\partial x_2}\left(\frac{1}{\sigma}\frac{\partial}{\partial x_2}
		\left(\frac{1}{\sigma}\frac{\partial\tilde{u}}{\partial x_2}\right)+\mu\frac{\partial^2\tilde{u}}{\partial x_1^2}\right)\right){\rm d}s
		+\int_{\Gamma_1}\frac{1}{\sigma} \frac{\partial \overline{v}}{\partial x_2}
	\left(\frac{1}{\sigma}\frac{\partial}{\partial x_2}
		\frac{1}{\sigma}\frac{\partial\tilde{u}}{\partial x_2}+\mu\frac{\partial^2\tilde{u}}{\partial x_1^2}\right){\rm d}s\\
&\quad+\int_{\Omega^{\rm PML}_1}\frac{\partial}{\partial x_2}\left( \frac{1}{\sigma} \frac{\partial \overline{v}}{\partial x_2}\right)\,
		\left(\frac{1}{\sigma}\frac{\partial}{\partial x_2}
		\left(\frac{1}{\sigma}\frac{\partial\tilde{u}}{\partial x_2}\right)+\mu\frac{\partial^2\tilde{u}}{\partial x_1^2}\right){\rm d}x.
\end{align*}
Combining the above results leads to 
\begin{align*}
& \int_{\Omega^{\rm PML}_1} \Bigg[\frac{\partial^2}{\partial x^2_1}\left(\sigma\frac{\partial^2 \tilde{u}}{\partial x_1^2}
	+\mu\frac{\partial}{\partial x_2}\left(\frac{1}{\sigma}\frac{\partial \tilde{u}}{\partial x_2}\right)\right)
	+2(1-\mu)\frac{\partial^2}{\partial x_1 \partial x_2}\left(\frac{1}{\sigma}\frac{\partial^2 \tilde{u}}{\partial x_1 \partial x_2}\right)\\
&\quad	+\frac{\partial}{\partial x_2}\left(
	\frac{1}{\sigma}\frac{\partial}{\partial x_2}\left(\frac{1}{\sigma}\frac{\partial}{\partial x_2}
		\left(\frac{1}{\sigma}\frac{\partial \tilde{u}}{\partial x_2}\right)+\mu\frac{\partial^2 \tilde{u}}{\partial x_1^2}\right)\right)	
		-\kappa^4 \sigma \tilde{u}\Bigg]\overline{v}{\rm d}x\\
&=	\int_{\Omega^{\rm PML}_1} \Bigg[\frac{\partial^2 \overline{v}}{\partial x_1^2}\,\left(\sigma\frac{\partial^2 \tilde{u}}{\partial x_1^2}
	+\mu\frac{\partial}{\partial x_2}\left(\frac{1}{\sigma}\frac{\partial\tilde{u}}{\partial x_2}\right)\right)
	+\frac{2-2\mu}{\sigma}\frac{\partial^2\overline{v}}{\partial x_1 \partial x_2}\,\frac{\partial^2 \tilde{u}}{\partial x_1 \partial x_2}\\
&\quad	+\frac{\partial}{\partial x_2}\left( \frac{1}{\sigma} \frac{\partial \overline{v}}{\partial x_2}\right)\,
		\left(\frac{1}{\sigma}\frac{\partial}{\partial x_2}
		\left(\frac{1}{\sigma}\frac{\partial\tilde{u}}{\partial x_2}\right)
		+\mu\frac{\partial^2\tilde{u}}{\partial x_1^2}\right)\Bigg]{\rm d}x
	-\int_{\Gamma_1}(2-2\mu)\overline{v}\frac{\partial^3 \tilde{u}}{\partial x_1^2 \partial x_2}\,{\rm d}s	\\
&\quad-\int_{\Gamma_1} \overline{v}\,	
	\frac{\partial}{\partial x_2}\left(\frac{1}{\sigma}\frac{\partial}{\partial x_2}
		\left(\frac{1}{\sigma}\frac{\partial\tilde{u}}{\partial x_2}\right)+\mu\frac{\partial^2\tilde{u}}{\partial x_1^2}\right){\rm d}s
		+\int_{\Gamma_1} \frac{\partial \overline{v}}{\partial x_2}
	\left(\frac{\partial}{\partial x_2}\left(
		\frac{1}{\sigma}\frac{\partial\tilde{u}}{\partial x_2}\right)+\mu\frac{\partial^2\tilde{u}}{\partial x_1^2}\right){\rm d}s. 
\end{align*}

Since $\sigma$ is assumed to be three times differentiable, we have
\begin{align*}
&-\int_{\Gamma_1} \overline{v}\,	
	\frac{\partial}{\partial x_2}\left(\frac{1}{\sigma}\frac{\partial}{\partial x_2}
		\left(\frac{1}{\sigma}\frac{\partial\tilde{u}}{\partial x_2}\right)+\mu\frac{\partial^2\tilde{u}}{\partial x_1^2}\right){\rm d}s
		+\int_{\Gamma_1} \frac{\partial \overline{v}}{\partial x_2}
	\left(\frac{\partial}{\partial x_2}\left(
		\frac{1}{\sigma}\frac{\partial\tilde{u}}{\partial x_2}\right)+\mu\frac{\partial^2\tilde{u}}{\partial x_1^2}\right){\rm d}s\\
&=\int_{\Gamma_1}\Bigg[ -\mu\,\overline{v}\,\frac{\partial^3 \tilde{u}}{\partial x_1^2 \partial x_2}
	-\left(\frac{\partial^3 \tilde{u}}{\partial x_2^3}-3\sigma'\frac{\partial^2 \tilde{u}}{\partial x_2^2}
	-\left(\sigma''-3\sigma'^2\right)\frac{\partial \tilde{u}}{\partial x_2}\right)\overline{v}
	+\left(-\sigma'\frac{\partial \tilde{u}}{\partial x_2}+\frac{\partial^2 \tilde{u}}{\partial x_2^2}
		+\mu\frac{\partial^2 \tilde{u}}{\partial x_1^2}\right) \frac{\partial \overline{v}}{\partial x_2}\Bigg]{\rm d}s.
\end{align*}
Considering that $\sigma'(h_1)=\sigma''(h_1)=0$, it holds that
\begin{align}\label{Upper}
 &\int_{\Omega^{\rm PML}_1} \Bigg[\frac{\partial^2}{\partial x^2_1}\left(\sigma\frac{\partial^2 \tilde{u}}{\partial x_1^2}
	+\mu\frac{\partial}{\partial x_2}\left(\frac{1}{\sigma}\frac{\partial \tilde{u}}{\partial x_2}\right)\right)
	+2(1-\mu)\frac{\partial^2}{\partial x_1 \partial x_2}
	\left(\frac{1}{\sigma}\frac{\partial^2 \tilde{u}}{\partial x_1 \partial x_2}\right)\notag\\
&\quad +\frac{\partial}{\partial x_2}\left(
	\frac{1}{\sigma}\frac{\partial}{\partial x_2}\left(\frac{1}{\sigma}\frac{\partial}{\partial x_2}
		\left(\frac{1}{\sigma}\frac{\partial \tilde{u}}{\partial x_2}\right)+\mu\frac{\partial^2 \tilde{u}}{\partial x_1^2}\right)\right)	
		-\kappa^4 \sigma \tilde{u}\Bigg]\overline{v}{\rm d}x\notag\\
&=	\int_{\Omega^{\rm PML}_1} \Bigg[\frac{\partial^2 \overline{v}}{\partial x_1^2}\,\left(\sigma\frac{\partial^2 \tilde{u}}{\partial x_1^2}
	+\mu\frac{\partial}{\partial x_2}\left(\frac{1}{\sigma}\frac{\partial\tilde{u}}{\partial x_2}\right)\right)
	+\frac{2-2\mu}{\sigma}\frac{\partial^2\overline{v}}{\partial x_1 \partial x_2}\,\frac{\partial^2 \tilde{u}}{\partial x_1 \partial x_2}\notag\\
&\quad +\frac{\partial}{\partial x_2}\left( \frac{1}{\sigma} \frac{\partial \overline{v}}{\partial x_2}\right)\,
		\left(\frac{1}{\sigma}\frac{\partial}{\partial x_2}
		\left(\frac{1}{\sigma}\frac{\partial\tilde{u}}{\partial x_2}\right)
		+\mu\frac{\partial^2\tilde{u}}{\partial x_1^2}\right)-\kappa^4 \sigma \tilde{u}\,\overline{v}\Bigg]{\rm d}x\notag\\
&\quad +\int_{\Gamma_1}\left(\overline{v} N_1 \tilde{u}+\partial_{\nu} \overline{v} M_1 \tilde{u}\right){\rm d}s.
\end{align}

Following the same argument, we deduce that 
\begin{align}\label{Lower}
& \int_{\Omega^{\rm PML}_2} \Bigg[\frac{\partial^2}{\partial x^2_1}\left(\sigma\frac{\partial^2 \tilde{u}}{\partial x_1^2}+\mu\frac{\partial}{\partial x_2}\left(\frac{1}{\sigma}\frac{\partial \tilde{u}}{\partial x_2}\right)\right)
	+2(1-\mu)\frac{\partial^2}{\partial x_1 \partial x_2}
		\left(\frac{1}{\sigma}\frac{\partial^2 \tilde{u}}{\partial x_1 \partial x_2}\right)\notag\\
&\quad +\frac{\partial}{\partial x_2}\left(
	\frac{1}{\sigma}\frac{\partial}{\partial x_2}\left(\frac{1}{\sigma}\frac{\partial}{\partial x_2}
		\left(\frac{1}{\sigma}\frac{\partial \tilde{u}}{\partial x_2}\right)+\mu\frac{\partial^2 \tilde{u}}{\partial x_1^2}\right)\right)
		-\kappa^4 \sigma \tilde{u}\Bigg]\overline{v}{\rm d}x\notag\\
& =\int_{\Omega^{\rm PML}_2} 
	\Bigg[\frac{\partial^2 \overline{v}}{\partial x_1^2}\,\left(\sigma\frac{\partial^2 \tilde{u}}{\partial x_1^2}
	+\mu\frac{\partial}{\partial x_2}\left(\frac{1}{\sigma}\frac{\partial\tilde{u}}{\partial x_2}\right)\right)
	+\frac{2-2\mu}{\sigma}\frac{\partial^2\overline{v}}{\partial x_1 \partial x_2}\,\frac{\partial^2 \tilde{u}}{\partial x_1 \partial x_2}\notag\\
&\quad  +\frac{\partial}{\partial x_2}\left( \frac{1}{\sigma} \frac{\partial \overline{v}}{\partial x_2}\right)\,
		\left(\frac{1}{\sigma}\frac{\partial}{\partial x_2}
		\left(\frac{1}{\sigma}\frac{\partial\tilde{u}}{\partial x_2}\right)
		+\mu\frac{\partial^2\tilde{u}}{\partial x_1^2}\right)-\kappa^4 \sigma \tilde{u}\,\overline{v}\Bigg]{\rm d}x\notag\\
&\quad +\int_{\Gamma_2}\left(\overline{v} N_2 \tilde{u}+\partial_{\nu} \overline{v} M_2 \tilde{u}\right){\rm d}s.
\end{align}

Substituting \eqref{Physical}--\eqref{Lower} into \eqref{auv}, we arrive at the variational problem for \eqref{MBHE}: to find $\tilde{u}\in H^2_{{\rm qp}, 0}(\Omega^{\rm PML})$ with $\tilde{u}=u^i, \partial_{\nu}\tilde{u}=\partial_{\nu} u^i$ on $\Gamma_1^{\rm PML}$ such that
\begin{equation}\label{variationalP1}
a^{\rm PML}(\tilde{u}, v)=\int_{\Omega^{\rm PML}} \sigma f\overline{v}{\rm d}x,\quad\forall\, v\in H^2_{{\rm qp}, 00}(\Omega^{\rm PML}), 
\end{equation}
where the sesquilinear form $a^{\rm PML}(u, v): H_{{\rm qp}, \Gamma_c}^2(\Omega^{\rm PML})\times
H_{{\rm qp}, \Gamma_c}^2(\Omega^{\rm PML})\rightarrow\mathbb{C}$ is defined as
\begin{align}\label{auvpml}
a^{\rm PML}(\tilde{u}, v) &=
\int_{\Omega^{\rm PML}} \Bigg[\sigma\frac{\partial^2 \tilde{u}}{\partial x_1^2}\frac{\partial^2 \overline{v}}{\partial x_1^2}+\mu\frac{\partial}{\partial x_2}\left(\frac{1}{\sigma}\frac{\partial\tilde{u}}{\partial x_2}\right)\frac{\partial^2 \overline{v}}{\partial x_1^2}+\frac{2-2\mu}{\sigma}\frac{\partial^2\overline{v}}{\partial x_1 \partial x_2}\frac{\partial^2 \tilde{u}}{\partial x_1 \partial x_2}\notag\\
&\quad +\frac{1}{\sigma}\frac{\partial}{\partial x_2}\left(\frac{1}{\sigma}\frac{\partial\tilde{u}}{\partial x_2}\right)
		\frac{\partial}{\partial x_2}\left( \frac{1}{\sigma} \frac{\partial \overline{v}}{\partial x_2}\right)
		+\mu\frac{\partial^2\tilde{u}}{\partial x_1^2}\frac{\partial}{\partial x_2}
		\left( \frac{1}{\sigma} \frac{\partial \overline{v}}{\partial x_2}\right)-\kappa^4 \sigma \tilde{u}\overline{v}\Bigg]{\rm d}x.
\end{align}

\subsection{The well-posedness}\label{PML:wellp}

The remainder of this section aims to establish the well-posedness of the variational problem \eqref{variationalP1}. Let
$\sigma(t)=\alpha_1(t)+{\rm i}\alpha_2(t).$ By letting $v=\tilde{u}$ and taking the real part of \eqref{auvpml}, we get
\begin{eqnarray*}
\Re\int_{\Omega^{\rm PML}} 
	\sigma\frac{\partial^2 \tilde{u}}{\partial x_1^2}\overline{\frac{\partial^2 \tilde{u}}{\partial x_1^2}}
	+\frac{2-2\mu}{\sigma}\overline{\frac{\partial^2\tilde{u}}{\partial x_1 \partial x_2}}\,
	\frac{\partial^2 \tilde{u}}{\partial x_1 \partial x_2}{\rm d}x
=\int_{\Omega^{\rm PML}} \alpha_1 \left|\frac{\partial^2 \tilde{u}}{\partial x_1^2}\right|^2
+\frac{\alpha_1(2-2\mu)}{|\sigma|^2}\left|\frac{\partial^2\tilde{u}}{\partial x_1 \partial x_2}\right|	^2{\rm d}x,
\end{eqnarray*}
and
\begin{align*}
&\Re\int_{\Omega^{\rm PML}} \left[
\mu\frac{\partial}{\partial x_2}\left(\frac{1}{\sigma}\frac{\partial\tilde{u}}{\partial x_2}\right)\,
	\overline{\frac{\partial^2 \tilde{u}}{\partial x_1^2}}
	+\mu\frac{\partial^2\tilde{u}}{\partial x_1^2}\frac{\partial}{\partial x_2}
		\left( \frac{1}{\sigma} \overline{\frac{\partial \tilde{u}}{\partial x_2}}\right)
		+\frac{1}{\sigma}\frac{\partial}{\partial x_2}
		\left(\frac{1}{\sigma}\frac{\partial\tilde{u}}{\partial x_2}\right)
		\frac{\partial}{\partial x_2}\left( \frac{1}{\sigma}\overline{ \frac{\partial \tilde{u}}{\partial x_2}}\right)\right]{\rm d}x\\
&=\Re\int_{\Omega^{\rm PML}} \Bigg[-\mu\frac{\sigma'}{\sigma^2}
	\left(\frac{\partial\tilde{u}}{\partial x_2}\,
	\overline{\frac{\partial^2 \tilde{u}}{\partial x_1^2}}+\overline{\frac{\partial\tilde{u}}{\partial x_2}}\,
	\frac{\partial^2 \tilde{u}}{\partial x_1^2}\right)
	+\frac{\mu}{\sigma}\left(\frac{\partial^2 \tilde{u}}{\partial x_1^2}\overline{\frac{\partial^2 \tilde{u}}{\partial x_2^2}}
	+\overline{\frac{\partial^2 \tilde{u}}{\partial x_1^2}}\frac{\partial^2 \tilde{u}}{\partial x_2^2}\right)\\
&\quad+\frac{1}{\sigma}\left(-\frac{\sigma'}{\sigma^2}\frac{\partial\tilde{u}}{\partial x_2}
	+\frac{1}{\sigma}\frac{\partial^2\tilde{u}}{\partial x^2_2}\right)
\left(-\frac{\sigma'}{\sigma^2}\overline{\frac{\partial\tilde{u}}{\partial x_2}}
+\frac{1}{\sigma}\overline{\frac{\partial^2\tilde{u}}{\partial x^2_2}}\right)
	\Bigg]{\rm d}x\\
&=\int_{\Omega^{\rm PML}}\Bigg[
	-\mu\,\Re\frac{\sigma'}{\sigma^2}
	\left(\frac{\partial\tilde{u}}{\partial x_2}\,
	\overline{\frac{\partial^2 \tilde{u}}{\partial x_1^2}}+\overline{\frac{\partial\tilde{u}}{\partial x_2}}\,
	\frac{\partial^2 \tilde{u}}{\partial x_1^2}\right)
	+\frac{\mu\alpha_1}{|\sigma|^2}
	\left(\frac{\partial^2 \tilde{u}}{\partial x_1^2}\overline{\frac{\partial^2 \tilde{u}}{\partial x_2^2}}
	+\overline{\frac{\partial^2 \tilde{u}}{\partial x_1^2}}\frac{\partial^2 \tilde{u}}{\partial x_2^2}\right)\\
&\quad
	+\Re\frac{1}{\sigma^3}\left|\frac{\partial^2\tilde{u}}{\partial x^2_2}\right|^2
	+\Re\frac{\sigma'^2}{\sigma^5}\left|\frac{\partial\tilde{u}}{\partial x_2}\right|^2	
	-\Re\frac{\sigma'}{\sigma^4}\left(\frac{\partial\tilde{u}}{\partial x_2}\overline{\frac{\partial^2\tilde{u}}{\partial x^2_2}}
	+\overline{\frac{\partial\tilde{u}}{\partial x_2}}\frac{\partial^2\tilde{u}}{\partial x^2_2}\right)
\Bigg]{\rm d}x.
\end{align*}

Taking any positive values $\epsilon_1, \epsilon_2$, and $\epsilon_3$, we have from Cauchy's inequality that
\begin{align*}
\Re a^{\rm PML}(\tilde{u}, \tilde{u})\geq \int_{\Omega^{\rm PML}}\Bigg[\left( \alpha_1-\epsilon_1-\frac{\mu\alpha_1 \epsilon_2}{|\sigma|^2}\right)\left|\frac{\partial^2 \tilde{u}}{\partial x_1^2}\right|^2
+(2-2\mu)\frac{\alpha_1}{|\sigma|^2}\left|\frac{\partial^2\tilde{u}}{\partial x_1 \partial x_2}\right|^2\\
+\left(\Re\frac{1}{\sigma^3}-\frac{\mu\alpha_1}{|\sigma|^2 \epsilon_2}-\epsilon_3\right) 
\left|\frac{\partial^2\tilde{u}}{\partial x^2_2}\right|^2\Bigg]{\rm d}x
-\int_{\Omega^{\rm PML}}\Bigg(\left|-\mu \Re\frac{\sigma'}{\sigma^2}\right|^2\frac{1}{\epsilon_1}
+\left|\Re\frac{\sigma'^2}{\sigma^5}\right|\\
+\left|\Re\frac{\sigma'}{\sigma^4}\right|^2\frac{1}{\epsilon_3}
\Bigg)\left|\frac{\partial \tilde{u}}{\partial x_2}\right|^2{\rm d}x
-\int_{\Omega^{\rm PML}} \kappa^4 \alpha_1 |\tilde{u}|^2{\rm d}x.
\end{align*}

Taking $\mu=\frac{1}{2}$ and $\sigma_2=\frac{1}{3}(1+\sigma_1)$ for clarity, we have 
\[
0<\frac{\alpha_2(t)}{\alpha_1(t)}=\frac{\sigma_2 t^m}{1+\sigma_1 t^m}<\frac{1}{3},
\quad 1< |\sigma(t)|^2<\frac{10}{9}(1+\sigma_1)^2,\quad \forall\, t\in[0, 1].
\]
Choosing $\epsilon_2=\frac{1}{2}\left(\frac{|\sigma|^2}{\mu}+\mu\frac{|\sigma|^4}{\alpha_1^2-3\alpha_2^2}\right)$ yields
\begin{align*}
\alpha_1-\frac{\mu\alpha_1 \epsilon_2}{|\sigma|^2}&=
\alpha_1\left(1-\frac{\mu}{2|\sigma|^2}\left(\frac{|\sigma|^2}{\mu}+\frac{\mu|\sigma|^4}{\alpha_1^2-3\alpha_2^2}\right)\right)
=\frac{\alpha_1}{2}\left(1-\mu^2-\frac{4\alpha_2^2 \mu^2}{\alpha_1^2-3\alpha_2^2}\right)\\
&>\frac{1}{2}\left(1-\frac{1}{4}-\frac{4\left(\frac{\alpha_2}{\alpha_1}\right)^2 \mu^2}{1-3\left(\frac{\alpha_2}{\alpha_1}\right)^2}\right)
>\frac{7}{24},
\end{align*}
and
\begin{align*}
\Re\frac{1}{\sigma^3}-\frac{\mu\alpha_1}{|\sigma|^2 \epsilon_2}&=
\frac{\alpha_1^3}{|\sigma|^6}-\frac{\mu\alpha_1}{|\sigma|^2}\frac{2}{\frac{|\sigma|^2}{\mu}+\frac{\mu|\sigma|^4}{\alpha_1^2-3\alpha_2^2}}
=\frac{\alpha_1(\alpha_1^2-3\alpha_2^2)}{|\sigma|^6}\left[1-\frac{2\mu^2}{\frac{\alpha_1^2-3\alpha_2^2}{\alpha_1^2+\alpha_2^2}+\mu^2}\right]\\
&>\frac{1}{\alpha_1^3}\frac{1-3\left(\frac{\alpha_2}{\alpha_1}\right)^2}{\left[1+\left(\frac{\alpha_2}{\alpha_1}\right)^2\right]^3}
\left[1-\frac{2\mu^2}{\frac{\alpha_1^2-3\alpha_2^2}{\alpha_1^2+\alpha_2^2}+\mu^2}\right]
>\frac{14}{51}\left(\frac{9}{10}\right)^3\left(\frac{1}{1+\sigma_1}\right)^3.
\end{align*}
Taking the parameter $\epsilon_1$ and $\epsilon_2$  as
\[
	\epsilon_1=\frac{1}{2}\left(\alpha_1-\frac{\mu\alpha_1 \epsilon_2}{|\sigma|^2}\right)>\frac{7}{48},\quad
	\epsilon_3=\frac{1}{2}\left(\Re\frac{1}{\sigma^3}-\frac{\mu\alpha_1}{|\sigma|^2 \epsilon_2}\right)
	>\frac{7}{51}\left(\frac{9}{10}\right)^3\left(\frac{1}{1+\sigma_1}\right)^3,
\]
we obtain that the coefficient of $\partial_{x_2}\tilde{u}$ is bounded by
\begin{eqnarray*}
\left|-\mu \Re\frac{\sigma'}{\sigma^2}\right|^2\frac{1}{\epsilon_1}
+\left|\Re\frac{\sigma'^2}{\sigma^5}\right|
+\left|\Re\frac{\sigma'}{\sigma^4}\right|^2\frac{1}{\epsilon_3}\leq
\frac{24}{7}\frac{|\sigma'|^2}{|\sigma^2|^2}
+\frac{|\sigma'^2|}{|\sigma^5|}+\frac{51}{7}\left(\frac{10}{9}\right)^3\frac{(1+\sigma_1)^3|\sigma'|^2}{|\sigma^4|^2}\\
\leq \left(\frac{31}{7}+\frac{51(1+\sigma_1)^3}{7}\left(\frac{10}{9}\right)^3\right)m^2(\sigma_1^2+\sigma_2^2).
\end{eqnarray*}
It follows from the above estimates that the sesquilinear form \eqref{auvpml} satisfies G\r{a}rding's inequality
\[
	\Re a(\tilde{u}, \tilde{u})\geq 
	C_1|\tilde u|^2_{H^2(\Omega^{\rm PML})}
	-C_2 \|\tilde u\|^2_{H^1(\Omega^{\rm PML})},
\]
where
\begin{eqnarray*}
 &&C_1=\min \left\{\frac{7}{24}, \,\frac{7}{51}\left(\frac{9}{10}\right)^3\left(\frac{1}{1+\sigma_1}\right)^3,\,
 \frac{9}{20}\frac{1}{(1+\sigma_1)^2}\right\},\\
 &&C_2=\max\left\{ \left(\frac{31}{7}+\frac{51(1+\sigma_1)^3}{7}\left(\frac{10}{9}\right)^3\right)m^2(\sigma_1^2+\sigma_2^2),\,
 \kappa^2 (1+\sigma_1)\right\}.
\end{eqnarray*}
The well-posedness of the variational problem \eqref{variationalP1} can be proven by applying the Fredholm alternative theorem, as summarized in the following result.

\begin{theorem}\label{PMLwellposed}
For any $\mu\in(0, 1)$, if the PML parameter $\sigma$ is chosen according to \eqref{pmlcons} with $m>3$ and satisfies \eqref{PMLrestriction}, then the variational problem \eqref{variationalP1} admits a unique weak solution except for a discrete set of wavenumbers $\kappa$.
\end{theorem}

\section{Convergence analysis}\label{section5}

This section focuses on an error estimate of the solutions to the truncated PML problem \eqref{variationalP1} and the original scattering problem \eqref{TotalBiharmonicReduced}.

Introduce a reference function $\hat{u}\in H^2_{\rm qp}(\Omega^{\rm PML})$, which corresponds to the solution $u$ of the
problem \eqref{TotalBiharmonicReduced} in $\Omega$. This function has the Fourier series expansion in the upper PML layer $\Omega_1^{\rm PML}$:
\begin{align}\label{hatu1}
\hat{u}(x_1, x_2)&=u^i(x_1, x_2)+\sum\limits_{n\in\mathbb{Z}} \left(
\frac{\gamma_n \hat{f}^{(n)}_1+\hat{g}^{(n)}_1}{\gamma_n+{\rm i}\beta_n}\right)e^{{\rm i}\alpha_n x_1
+{\rm i}\beta_n(\tilde{x}_2(x_2)-h_1)}\notag\\
&\quad +\sum\limits_{n\in\mathbb{Z}}\left( \frac{{\rm i}\beta_n \hat{f}^{(n)}_1-\hat{g}^{(n)}_1}{\gamma_n+{\rm i}\beta_n}\right)e^{{\rm i}\alpha_n x_1-\gamma_n(\tilde{x}_2(x_2)-h_1)},
\end{align}
and the Fourier series expansion in the lower PML layer $\Omega_2^{\rm PML}$:
\begin{eqnarray}
\hat{u}(x_1, x_2)=\sum\limits_{n\in\mathbb{Z}} \left[\left(
		\frac{\gamma_n f^{(n)}_2+g^{(n)}_2}{\gamma_n+{\rm i}\beta_n}\right)
			e^{-{\rm i}\beta_n(\tilde{x}_2(x_2)-h_2)}
			+\left(\frac{{\rm i}\beta_n f^{(n)}_2-g^{(n)}_2}{\gamma_n+{\rm i}\beta_n}
		\right)e^{\gamma_n(\tilde{x}_2(x_2)-h_2)}\right] e^{{\rm i}\alpha_n x_1},\label{hatu2}
\end{eqnarray}
where $\big(\hat{f}_1^{(n)}, \hat{g}^{(n)}_1\big)$ and $\big(f_2^{(n)}, g^{(n)}_2\big)$ represent the Fourier coefficients of the Dirichlet data $\left(u-u^i, \partial_{\nu} u-\partial_{\nu} u^i\right)$ and
 $\left(u, \partial_\nu u\right)$ on $\Gamma_1$ and $\Gamma_2$, respectively.

It can be verified that $\hat{u}$ satisfies the following system
\begin{equation}\label{hatu3}
\left\{
\begin{aligned}
&\tilde{\Delta}^2\hat{u}-\kappa^4\hat{u}=f & &{\rm in} ~ \Omega^{\rm PML},\\
&\hat{u}=P_1(\hat{f}_1, \hat{g}_1)+u^i  & &{\rm on} ~ \Gamma_1^{\rm PML},\\
&\partial_{\nu}\hat{u}=Q_1(\hat{f}_1, \hat{g}_1)+\partial_{\nu} u^i & &{\rm on} ~ \Gamma_1^{\rm PML},\\
&\hat{u}=P_2(f_2, g_2),\quad\partial_{\nu}\hat{u}=Q_2(f_2, g_2) & &{\rm on} ~ \Gamma_2^{\rm PML},\\
&\hat{u}=0,\quad\partial_{\nu}\hat{u}=0&&{\rm on}~\Gamma_c,
\end{aligned}
\right.
\end{equation}
where the propagating operators $P_i, Q_i$ are defined by
\begin{align}
P_1(\hat{f}_1, \hat{g}_1) & =
\sum\limits_{n\in\mathbb{Z}} \left(
		\frac{\gamma_n \hat{f}^{(n)}_1+\hat{g}^{(n)}_1}{\gamma_n+{\rm i}\beta_n}\right)e^{{\rm i}\alpha_n x_1
		+{\rm i}\beta_n(1+\frac{\sigma_1}{m+1}+{\rm i}\,\frac{\sigma_2}{m+1})\Delta h_1}\label{P1}\\
&\quad		+\sum\limits_{n\in\mathbb{Z}}
			\left( \frac{{\rm i}\beta_n \hat{f}^{(n)}_1-\hat{g}^{(n)}_1}{\gamma_n+{\rm i}\beta_n}\right)
				e^{{\rm i}\alpha_n x_1-\gamma_n(1+\frac{\sigma_1}{m+1}+{\rm i}\,\frac{\sigma_2}{m+1})\Delta h_1}\notag,\\
Q_1(\hat{f}_1, \hat{g}_1) & =
\sum\limits_{n\in\mathbb{Z}}{\rm i}\,\beta_n (1+\sigma_1+{\rm i}\,\sigma_2) \left(
		\frac{\gamma_n \hat{f}^{(n)}_1+\hat{g}^{(n)}_1}{\gamma_n+{\rm i}\beta_n}\right)e^{{\rm i}\alpha_n x_1
		+{\rm i}\beta_n(1+\frac{\sigma_1}{m+1}+{\rm i}\,\frac{\sigma_2}{m+1})\Delta h_1}\label{Q1}\\
&\quad	-\sum\limits_{n\in\mathbb{Z}}\gamma_n (1+\sigma_1+{\rm i}\,\sigma_2)
			\left( \frac{{\rm i}\beta_n \hat{f}^{(n)}_1-\hat{g}^{(n)}_1}{\gamma_n+{\rm i}\beta_n}\right)
				e^{{\rm i}\alpha_n x_1-\gamma_n(1+\frac{\sigma_1}{m+1}+{\rm i}\,\frac{\sigma_2}{m+1})\Delta h_1}\notag,
\end{align}
and
\begin{align}	
P_2(f_2, g_2) &=\sum\limits_{n\in\mathbb{Z}} \left(
		\frac{\gamma_n f^{(n)}_2+g^{(n)}_2}{\gamma_n+{\rm i}\beta_n}\right)
			e^{{\rm i}\alpha_n x_1-{\rm i}\beta_n(1+\frac{\sigma_1}{m+1}+{\rm i}\,\frac{\sigma_2}{m+1})\Delta h_2}\label{P2}\\
&\quad	+\sum\limits_{n\in\mathbb{Z}} \left(\frac{{\rm i}\beta_n f^{(n)}_2-g^{(n)}_2}{\gamma_n+{\rm i}\beta_n}
		\right)e^{{\rm i}\alpha_n x_1+\gamma_n(1+\frac{\sigma_1}{m+1}+{\rm i}\,\frac{\sigma_2}{m+1})\Delta h_2},\notag	\\
Q_2(f_2, g_2) & =\sum\limits_{n\in\mathbb{Z}} {\rm i}\,\beta_n (1+\sigma_1+{\rm i}\,\sigma_2)\left(
		\frac{\gamma_n f^{(n)}_2+g^{(n)}_2}{\gamma_n+{\rm i}\beta_n}\right)
			e^{{\rm i}\alpha_n x_1-{\rm i}\beta_n(1+\frac{\sigma_1}{m+1}+{\rm i}\,\frac{\sigma_2}{m+1})\Delta h_2}\label{Q2}\\
&\quad-\sum\limits_{n\in\mathbb{Z}} \gamma_n (1+\sigma_1+{\rm i}\,\sigma_2)
\left(\frac{{\rm i}\beta_n f^{(n)}_2-g^{(n)}_2}{\gamma_n+{\rm i}\beta_n}
		\right)e^{{\rm i}\alpha_n x_1+\gamma_n(1+\frac{\sigma_1}{m+1}+{\rm i}\,\frac{\sigma_2}{m+1})\Delta h_2}.\notag	
\end{align}
The properties of the propagating operators $P_i, Q_i$, $i=1, 2$, are discussed in Lemma \ref{Step1}.

When $m>3$, it can be readily verified from \eqref{hatu1}--\eqref{hatu2} that $\hat{u}$ satisfies the same TBCs \eqref{TotalBiharmonicReduced} on $\Gamma_k, k=1, 2$: 
\begin{equation}\label{hatuTBC}
\left\{
\begin{aligned}
	& N_1 \hat{u}=T_{11}^{(1)} \hat{f}_1+T_{12}^{(1)} \hat{g}_1+p_1 &&\text{on} ~ \Gamma_1,\\
	& M_1 \hat{u}=T_{21}^{(1)} \hat{f}_1+T_{22}^{(1)} \hat{g}_1+p_2 &&\text{on} ~ \Gamma_1,\\
	& N_2 \hat{u}=T_{11}^{(2)} f_2+T_{12}^{(2)}g_2 & &\text{on} ~ \Gamma_2,\\
	& M_2 \hat{u}=T_{21}^{(2)}f_2+T_{22}^{(2)}g_2 & &\text{on} ~ \Gamma_2,
\end{aligned}
\right.
\end{equation}
where the DtN operators $T^{(k)}_{i,j}$ and the functions $f_k, g_k, p_k$ are defined in \eqref{T1}--\eqref{p12}. This enables us to transfer the error analysis from $\tilde{u}$ to $u$ to the error between $\tilde{u}$ and $\hat{u}$, both of which are defined in the same computational domain.

Since $\hat{u}=u$ in $\Omega$, we define $e=\hat{u}-\tilde{u}$ as the error between the truncated PML solution and the true solution. Let $\varphi_i=e|_{\Gamma_i}, \psi_i=\partial_{\nu} e|_{\Gamma_i}, \tilde{f}_i=\tilde{u}|_{\Gamma_i}, \tilde{g}_i=\partial_{\nu} \tilde{u}|_{\Gamma_i}, i=1,2,$ respectively. It is easy to check that the error function $e$ satisfies the system 
\begin{equation}\label{err}
\left\{
\begin{aligned}
\Delta^2 e-\kappa^4 e=0\quad &{\rm in}\,\Omega,\\
e=0,\,\partial_{\nu} e=0\quad &{\rm on}\,\Gamma_c,
\end{aligned}
\right.
\end{equation}
and the boundary condition on $\Gamma_1$:
\begin{align}\label{errGamma1}
&N_1 e-T_{11}^{(1)} \varphi_1-T_{12}^{(1)} \psi_1\notag\\
&=\left[N_1 \left(\hat{u}-u^i\right)-T_{11}^{(1)} \left(\hat{u}-u^i\right)
	-T_{12}^{(1)} \left(\partial_{\nu}\hat{u}-\partial_{\nu} u^i\right)\right]
+\left(N_1 u^i -T_{11}^{(1)} u^i-T_{12}^{(1)} \partial_{\nu} u^i \right)\notag\\
&\quad -\left(N_1 u^i-T_{11}^{(1)} u^i-T_{12}^{(1)} \partial_{\nu} u^i \right)
-\left[N_1 \left(\tilde{u}-u^i \right)-T_{11}^{(1)} \left(\tilde{u}-u^i\right)
	-T_{12}^{(1)} \left(\partial_{\nu}\tilde{u}-\partial_{\nu} u^{\rm i}\right)\right]\notag\\
& =-\left[N_1 \left(\tilde{u}-u^i \right)-T_{11}^{(1)} \left(\tilde{u}-u^{\rm i}\right)
	-T_{12}^{(1)} \left(\partial_{\nu}\tilde{u}-\partial_{\nu} u^i \right)\right].
\end{align}

Denote $f_1^*=\tilde{f}_1-u^i$ and $g_1^*=\tilde{g}_1-\partial_{\nu} u^i$. It follows from the continuity conditions across $\Gamma_1$ that
\[
	N_1 \left(\tilde{u}-u^i\right)=N_1 w_1^{(a)}, \quad
	T_{11}^{(1)} \left(\tilde{u}-u^i\right)
	+T_{12}^{(1)} \left(\partial_{\nu}\tilde{u}-\partial_{\nu} u^i \right)=N_1 w_1^{(b)},
\]
where $w^{(a)}_1$ and $w^{(b)}_1$ are the solutions of 
\[
\left\{
\begin{aligned}
&\tilde{\Delta}^2 w^{(a)}_1-\kappa^4 w^{(a)}_1=0 & &{\rm in} ~ \Omega_{1}^{\rm PML},\\
&w^{(a)}_1=f^*_1, \, \partial_{\nu}w^{(a)}_1=g^*_1& &{\rm on} ~ \Gamma_1,\\
&w^{(a)}_1=0 & &{\rm on} ~ \Gamma_1^{\rm PML},\\
&\partial_{\nu}w^{(a)}_1=0  & &{\rm on} ~ \Gamma_1^{\rm PML},\\
\end{aligned}
\right.
\]
and\[
\left\{
\begin{aligned}
&\tilde{\Delta}^2 w^{(b)}_1-\kappa^4 w^{(b)}_1=0 & &{\rm in} ~ \Omega_{1}^{\rm PML},\\
&w^{(b)}_1=f^*_1, \, \partial_{\nu}w^{(a)}_1=g^*_1& &{\rm on} ~ \Gamma_1,\\
&w^{(b)}_1=P_1(f^*_1, g^*_1) & &{\rm on} ~ \Gamma_1^{\rm PML},\\
&\partial_{\nu}w^{(b)}_1=Q_1(f^*_1, g^*_1)  & &{\rm on} ~ \Gamma_1^{\rm PML},\\
\end{aligned}
\right.
\]

Denoting  $w_{1}=w_1^{(b)}-w^{(a)}_1$, the boundary condition \eqref{errGamma1} is equivalent to
\begin{equation}\label{errGamma11}
	N_1 e-T_{11}^{(1)} \varphi_1-T_{12}^{(1)} \psi_1=N_1 w_1|_{\Gamma_1},
\end{equation}
where $w_1$ is the solution of the auxiliary problem
\begin{equation}\label{errGamma1b}
\left\{
\begin{aligned}
&\tilde{\Delta}^2 w_1-\kappa^4 w_1=0 & &{\rm in} ~ \Omega_{1}^{\rm PML},\\
&w_1=0, \, \partial_{\nu}w_1=0 & &{\rm on} ~ \Gamma_1,\\
&w_1=P_1(f^*_1, g^*_1), \quad \partial_{\nu}w_1=Q_1(f^*_1, g^*_1)  & &{\rm on} ~ \Gamma_1^{\rm PML}.
\end{aligned}
\right.
\end{equation}

Similarly, the TBCs for $e$ on $\Gamma_1$ and $\Gamma_2$ can be reformulated as
\begin{equation}\label{errGamma2}
\left\{
\begin{aligned}
	&M_1 e-T_{21}^{(1)} \varphi_1-T_{22}^{(1)} \psi_1 = M_1 w_1|_{\Gamma_1},\\
	&N_2 e-T_{11}^{(2)} \varphi_2-T_{12}^{(2)} \psi_2 = N_2 w_{2}|_{\Gamma_2},\\
	&M_2 e-T_{21}^{(2)} \varphi_2-T_{22}^{(2)} \psi_2 = M_2 w_{2}|_{\Gamma_2},
\end{aligned}
\right.
\end{equation}
where $w_2$ is the solution of the auxiliary problem
\begin{equation}\label{errGamma1c}
\left\{
\begin{aligned}
&\tilde{\Delta}^2 w_2-\kappa^4 w_2=0 & &{\rm in} ~ \Omega_{2}^{\rm PML},\\
&w_2=0, \quad \partial_{\nu} w_2=0& &{\rm on} ~ \Gamma_2,\\
&w_2=P_2(\tilde{f}_2, \tilde{g}_2),\quad \partial_{\nu}w_2=Q_2(\tilde{f}_2, \tilde{g}_2)  & &{\rm on} ~ \Gamma_2^{\rm PML}. 
\end{aligned}
\right.
\end{equation}

The equation \eqref{err}, combined with the TBCs \eqref{errGamma1} and \eqref{errGamma2}, leads to the variational problem for the error function: to find $e\in H^{2}_{{\rm qp}, \Gamma_c}(\Omega)$ such that 
\begin{align}\label{aev}
a(e, v)&=\int_{\Gamma_1} \left(\overline{v} N_1 w_1+\partial_{\nu}\overline{v} M_1 w_1\right){\rm d}x_1\notag\\
&\quad +\int_{\Gamma_2} \left(\overline{v} N_2 w_2+\partial_{\nu}\overline{v}_2 M_2 w_2\right){\rm d}x_1,\quad\forall\,v\in H^{2}_{{\rm qp}, \Gamma_c}(\Omega), 
\end{align}
where the sesquilinear form $a(\cdot, \cdot)$ is defined in \eqref{tbcauv}, and $w_k, k=1, 2,$ are the solutions of \eqref{errGamma1b} and \eqref{errGamma1c}, respectively. 

Inspired by the work presented in \cite{CL-SINUM-2005}, the framework for error analysis can be divided into three steps: (1) Lemma \ref{Step1} demonstrates that the propagating operators $P_i$ and $Q_i$ exponentially converge to zero with respect to the thickness of the layer. (2) Theorem \ref{Step2} establishes that the auxiliary problems \eqref{errGamma1b} and \eqref{errGamma1c} are well-defined, and shows that the boundary differential operators 
$M_i w_i$ and $N_i w_i$ can be bounded by the propagating operators. (3) Theorem \ref{Step3} verifies that variational problem \eqref{aev} is well-defined and that the error $e$ decays exponentially with respect to the thickness of the layer.

Define 
\[
\Delta^-=\min\limits_{n\in\mathbb{Z}}\left\{\Re\beta_n>0\right\}, \qquad
\Delta^+=\min\limits_{n\in\mathbb{Z}}\left\{\Im\beta_n>0\right\},
\]
and let $\delta=\min(\Delta h_1, \Delta h_2)$. Denote $\Theta$ as the PML efficiency constant given by
\begin{equation}\label{The}
\Theta=\max\left\{e^{-\frac{\sigma_2\delta}{m+1}\Delta^{-}}, e^{-\frac{m+1+\sigma_1}{m+1}\delta \Delta^{+}} ,e^{- \frac{m+1+\sigma_1}{m+1}\delta \sqrt{\kappa^2+\alpha^2}}\right\}.
\end{equation}

\begin{lemma}\label{Step1}
The propagating operators $P_k$, as defined in \eqref{P1} and \eqref{P2}, are bounded from 
$H^{3/2}(\Gamma_k)\times H^{1/2}(\Gamma_k)$ to $H^{3/2}(\Gamma_k^{\rm PML})$, and the propagating operators $Q_k$, defined in \eqref{Q1} and \eqref{Q2}, are bounded from $H^{3/2}(\Gamma_k)\times H^{1/2}(\Gamma_k)$ to $ H^{1/2}(\Gamma_k^{\rm PML})$, $k=1, 2$, respectively. Moreover, the following estimates hold 
\begin{align}
\|P_i(f_k, g_k)\|_{H^{3/2}\left(\Gamma_k^{\rm PML}\right)} &\lesssim \Theta \left(\|f_k\|_{H^{3/2}(\Gamma_k)}+\|g_k\|_{H^{1/2}(\Gamma_k)}\right),\label{ConvPi}\\
\|Q_k(f_k, g_k)\|_{H^{1/2}\left(\Gamma_k^{\rm PML}\right)} &\lesssim \Theta \left(\|f_k\|_{H^{3/2}(\Gamma_k)}+\|g_k\|_{H^{1/2}(\Gamma_k)}\right).\label{ConvQi}
\end{align}
\end{lemma}

\begin{proof}
It suffices to demonstrate the estimates for $P_1$ and $Q_1$, as the procedure is analogous to that of $P_2$ and $Q_2$. By \eqref{P1}, we have 
\begin{align*}
\|P_1(f_1, g_1)\|^2_{H^{3/2}\left(\Gamma_1^{\rm PML}\right)}&\leq \sum\limits_{n\in\mathbb{Z}}\left(1+n^2\right)^{3/2}
\left|\frac{\gamma_n f^{(n)}_1+g^{(n)}_1}{\gamma_n+{\rm i}\beta_n}\right|^2
\left|e^{{\rm i}\beta_n(1+\frac{\sigma_1}{m+1}+{\rm i}\,\frac{\sigma_2}{m+1})\Delta h_1}\right|^2\\
&\quad+ \sum\limits_{n\in\mathbb{Z}}\left(1+n^2\right)^{3/2}
\left| \frac{{\rm i}\beta_n f^{(n)}_1-g^{(n)}_1}{\gamma_n+{\rm i}\beta_n}\right|^2
\left|e^{-\gamma_n(1+\frac{\sigma_1}{m+1}+{\rm i}\,\frac{\sigma_2}{m+1})\Delta h_1}\right|^2\\
&\lesssim \Theta^2\sum\limits_{n\in\mathbb{Z}}\left[\left(1+n^2\right)^{3/2} |f_1^{(n)}|^2
+\left(1+n^2\right)^{1/2} |g_1^{(n)}|^2\right]\\
&\lesssim \Theta^2 \left(\|f_1\|^2_{H^{3/2}(\Gamma_1)}+\|g_1\|^2_{H^{1/2}(\Gamma_1)}\right),
\end{align*}
which verifies the estimate \eqref{ConvPi}. 

For the operator $Q_1$, a straightforward computation yields 
\begin{align*}
\|Q_1(f_1, g_1)\|^2_{H^{1/2}\left(\Gamma_1^{\rm PML}\right)}&\lesssim \sum\limits_{n\in\mathbb{Z}}\left(1+n^2\right)^{1/2}\Bigg[\left|{\rm i}\,\beta_n \left(\frac{\gamma_n \hat{f}^{(n)}_1+\hat{g}^{(n)}_1}{\gamma_n+{\rm i}\beta_n}\right)\right|^2 \left|e^{{\rm i}\beta_n(1+\frac{\sigma_1}{m+1}+{\rm i}\,\frac{\sigma_2}{m+1})\Delta h_1}\right|^2\\
&\quad+\left|\gamma_n\left( \frac{{\rm i}\beta_n \hat{f}^{(n)}_1-\hat{g}^{(n)}_1}{\gamma_n+{\rm i}\beta_n}\right)\right|^2\left|e^{-\gamma_n(1+\frac{\sigma_1}{m+1}+{\rm i}\,\frac{\sigma_2}{m+1})\Delta h_1}\right|^2\Bigg]\\
&\lesssim \Theta^2\sum\limits_{n\in\mathbb{Z}}\left[\left(1+n^2\right)^{3/2} |f_1^{(n)}|^2
+\left(1+n^2\right)^{1/2} |g_1^{(n)}|^2\right]\\
&\lesssim \Theta^2 \left(\|f_1\|^2_{H^{3/2}(\Gamma_1)}+\|g_1\|^2_{H^{1/2}(\Gamma_1)}\right),
\end{align*}
which confirms that the operator $Q_1$ is bounded from $H^{3/2}(\Gamma_1)\times H^{1/2}(\Gamma_1)$ to $H^{1/2}(\Gamma_1^{\rm PML})$ and converges to zero with respect to the thickness of the layer.
\end{proof}

\begin{theorem}\label{Step2}
For any given $f_k\in H^{3/2}(\Gamma_k)$ and $g_k\in H^{1/2}(\Gamma_k)$, if $\kappa$ is not an eigenvalue of the problems \eqref{errGamma1b} and \eqref{errGamma1c}, then the corresponding variational problems admit a unique weak solution in $H^{2}_{\rm qp}(\Omega_k^{\rm PML}), k=1, 2$. Furthermore, the following estimates hold:
\begin{equation}\label{estw}
\|w_k\|_{H^2\left(\Omega_k^{\rm PML}\right)}\lesssim \Theta \left(\|f_k\|_{H^{3/2}(\Gamma_k)}+\|g_k\|_{H^{1/2}(\Gamma_k)}\right),
\end{equation}
and
\begin{align}\label{estMN}
\|M_k w_k\|_{H^{-1/2}(\Gamma_k)} &\lesssim \Theta \left(\|f_k\|_{H^{3/2}(\Gamma_k)}+\|g_k\|_{H^{1/2}(\Gamma_k)}\right),\notag\\
\|N_k w_k\|_{H^{-3/2}(\Gamma_k)} &\lesssim \Theta \left(\|f_k\|_{H^{3/2}(\Gamma_k)}+\|g_k\|_{H^{1/2}(\Gamma_k)}\right).
\end{align}
\end{theorem}

\begin{proof}
Given $f_k\in H^{3/2}(\Gamma_k)$ and $g_k\in H^{1/2}(\Gamma_k)$, the variational problem for \eqref{errGamma1b}  and \eqref{errGamma1c} is to find $w_k\in H_{\rm qp}^2(\Omega_k^{\rm PML})$ satisfying $w_k=0, \partial_{\nu} w_k=0$ on $\Gamma_k$ and $w_k=P_k(f_k, g_k), \partial_{\nu} w_k=Q_k(p_k, g_k)$ on $\Gamma_k^{\rm PML}$ such that
\begin{equation}\label{buv}
b_k(w_k, v)=0, \quad\forall\, v\in H_{{\rm qp}, 00}^{2}(\Omega_k^{\rm PML}), 
\end{equation}
where the sesquilinear form $b_k: H_{\rm qp}^2(\Omega_k^{\rm PML})\times H_{\rm qp}^2(\Omega_k^{\rm PML})\rightarrow
\mathbb{C}$ is defined by
\begin{align*}
b_k(w, v)  &=\int_{\Omega_k^{\rm PML}} \Bigg[\sigma\frac{\partial^2 w }{\partial x_1^2}\frac{\partial^2 \overline{v}}{\partial x_1^2}+\mu\frac{\partial^2 \overline{v}}{\partial x_1^2}\frac{\partial}{\partial x_2}\left(\frac{1}{\sigma}\frac{\partial w}{\partial x_2}\right)+\frac{2-2\mu}{\sigma}\frac{\partial^2\overline{v}}{\partial x_1 \partial x_2}\frac{\partial^2 w}{\partial x_1 \partial x_2}\\
&\quad +\frac{1}{\sigma}\frac{\partial}{\partial x_2}\left( \frac{1}{\sigma} \frac{\partial \overline{v}}{\partial x_2}\right)\frac{\partial}{\partial x_2}\left(\frac{1}{\sigma}\frac{\partial w}{\partial x_2}\right)
+\mu\frac{\partial}{\partial x_2}\left( \frac{1}{\sigma} \frac{\partial \overline{v}}{\partial x_2}\right)\frac{\partial^2 w}{\partial x_1^2}-\kappa^4 \sigma w\overline{v}\Bigg]{\rm d}x.
\end{align*}
After establishing the proof of Theorem \ref{PMLwellposed}, it is evident that the sesquilinear form $b_k$ satisfies G\r{a}rding's inequality. Subsequently, utilizing the Fredholm alternative theorem, it can be shown that there exists a unique weak solution to  the variational problem \eqref{buv}, with the exception of a discrete set of wavenumbers $\kappa$.

Given that the PML parameter $|\sigma|$ has both a positive lower and upper bound in the PML regions, based on the general theory outlined in \cite[Chapter 5]{BA-AP-1973}, we deduce the existence of a constant $C>0$, such that
\begin{equation}\label{supb}
	\sup\limits_{0\neq\psi\in H_{{\rm qp}, 00}^{2}\left(\Omega^{\rm PML}_k\right)}
	\frac{|b_k(\varphi, \psi)|}{\|\psi\|_{H^2\left(\Omega^{\rm PML}_k\right)}}
	\geq C\|\varphi\|_{H^2\left(\Omega_k^{\rm PML}\right)},\quad
	\forall\varphi\in H_{{\rm qp}, 00}^{2}\left(\Omega_k^{\rm PML}\right).
\end{equation}

Drawing from the discussion provided in Appendix, we can identify a function 
$\varphi_k \in H_{\rm qp}^2(\Omega^{\rm PML}_k)$ satisfying boundary conditions such that $\varphi_k=0$ and 
$\partial_{\nu}\varphi_k=0$ on $\Gamma_k$, while $\varphi_k=P_k(f_k, g_k)$ and $\partial_{\nu}\varphi_k=Q_k(f_k, g_k)$ on $\Gamma^{\rm PML}_k$. Additionally, the following stability estimate holds: 
\[
\|\varphi_k\|_{H^2\left(\Omega_k^{\rm PML}\right)}\lesssim \|P_k(f_k, g_k)\|_{H^{3/2}\left(\Gamma_k^{\rm PML}\right)}
+\|Q_k(f_k, g_k)\|_{H^{1/2}\left(\Gamma_k^{\rm PML}\right)}.
\]

Let $\psi_k=w_k-\varphi_k\in H^2_{{\rm qp}, 00}\left(\Omega_k^{\rm PML}\right)$. Using the continuity of the sesquilinear form \eqref{buv} and the inf-sup condition \eqref{supb}, we deduce
\begin{align*}
\|\psi_k\|_{H^2\left(\Omega_k^{\rm PML}\right)} 
\lesssim \frac{|b_k(\psi_k, \psi_k)|}{\|\psi_k\|_{H^2(\Omega_k^{\rm PML})}}
= \frac{|b_k(\varphi_k, \psi_k)|}{\|\psi_k\|_{H^2(\Omega_k^{\rm PML})}}
\lesssim  \|\varphi_k\|_{H^2\left(\Omega_k^{\rm PML}\right)},
\end{align*}
which implies
\begin{align*}
	\|w_k\|_{H^2\left(\Omega_k^{\rm PML}\right)} &\leq 
	\|\psi_k\|_{H^2\left(\Omega_k^{\rm PML}\right)}+\|\varphi_k\|_{H^2\left(\Omega_k^{\rm PML}\right)}
	\lesssim  \|\varphi_k\|_{H^2\left(\Omega_k^{\rm PML}\right)}\\
	& \lesssim  \|P_k(f_k, g_k)\|_{H^{3/2}\left(\Gamma_k^{\rm PML}\right)}
	 +\|Q_k(f_k, g_k)\|_{H^{1/2}\left(\Gamma_k^{\rm PML}\right)}.
\end{align*}
Subsequently, the estimate \eqref{estw} can be derived by incorporating the estimates \eqref{ConvPi} and \eqref{ConvQi}.

To obtain the estimate \eqref{estMN}, following the discussion detailed in Appendix, we choose any $\varphi_k\in H^2_{\rm qp}\left(\Omega_k^{\rm PML}\right)$ satisfying the boundary conditions $\partial_\nu \varphi_k=0$ on $\Gamma_k$, and $\varphi_k=0, \partial_\nu \varphi_k=0$ on $\Gamma_k^{\rm PML}$ such that 
\[
	\|\varphi_k\|_{H^2\left(\Omega_k^{\rm PML}\right)}\lesssim \|\varphi_k\|_{H^{3/2}(\Gamma_k)}.
\]
Multiplying either \eqref{errGamma1b} or \eqref{errGamma1c} by $\varphi_k$ and applying integration by parts, we obtain 
\begin{align*}
\int_{\Gamma_k}\varphi_k N_k w_k{\rm d}x &=\int_{\Omega_k^{\rm PML}} \Bigg[
	\sigma\frac{\partial^2 w_k}{\partial x_1^2}\frac{\partial^2\varphi_k}{\partial x_1^2}
	+\mu\frac{\partial^2\varphi_k}{\partial x_1^2}\frac{\partial}{\partial x_2}\left(\frac{1}{\sigma}\frac{\partial w_k}{\partial x_2}\right)+\frac{2-2\mu}{\sigma}\frac{\partial^2\varphi_k}{\partial x_1 \partial x_2}\frac{\partial^2 w_k}{\partial x_1 \partial x_2}\\
&\quad +\frac{1}{\sigma}\frac{\partial}{\partial x_2}\left(\frac{1}{\sigma} \frac{\partial \varphi_k}{\partial x_2}\right)\frac{\partial}{\partial x_2}\left(\frac{1}{\sigma}\frac{\partial w_k}{\partial x_2}\right)
+\mu\frac{\partial}{\partial x_2}\left( \frac{1}{\sigma} \frac{\partial \varphi_k}{\partial x_2}\right)\frac{\partial^2 w_k}{\partial x_1^2}-\kappa^4 \sigma w_k\varphi_k\Bigg]{\rm d}x. 
\end{align*}
Thus
\begin{align*}
\left|\int_{\Gamma_k}\varphi_k N_k w_k{\rm d}x\right|&=\left|b_k(w_k, \varphi_k)\right|
\lesssim \|w_k\|_{H^2\left(\Omega_k^{\rm PML}\right)}\|\varphi_k\|_{H^2\left(\Omega_k^{\rm PML}\right)}\\
&\lesssim \|w_k\|_{H^2\left(\Omega_k^{\rm PML}\right)}\|\varphi_k\|_{H^{3/2}(\Gamma_k)}\\
&\lesssim \Theta \left(\|f_k\|_{H^{3/2}(\Gamma_k)}+\|g_k\|_{H^{1/2}(\Gamma_k)}\right)\|\varphi_k\|_{H^{3/2}(\Gamma_k)},
\end{align*}
which yields the estimate
\[
\|N_k w_k\|_{H^{-3/2}(\Gamma_k)}\lesssim\Theta \left(\|f_k\|_{H^{3/2}(\Gamma_k)}+\|g_k\|_{H^{1/2}(\Gamma_k)}\right).
\]

The estimate for $M_k$ can be obtained through a similar discussion by choosing any $\varphi_k\in H^2_{\rm qp}\left(\Omega_k^{\rm PML}\right)$ such that $\varphi_k=0$ on $\Gamma_k$, and $\varphi_k=0$, $\partial_\nu \varphi_k=0$ on $\Gamma_k^{\rm PML}$, and the following estimate holds:
\[
\|\varphi_k\|_{H^2(\Omega_k^{\rm PML})}\lesssim \|\partial_{\nu}\varphi_k\|_{H^{1/2}(\Gamma_k)}.
\]
which completes the proof. 
\end{proof}

The following result concerns the convergence of the truncated PML solution to that of the original scattering problem. 

\begin{theorem}\label{Step3}
If $e$ is the solution of \eqref{err} subject to the boundary conditions \eqref{errGamma11} and \eqref{errGamma2}, then
the variational problem  \eqref{aev} admits a unique solution in $H_{{\rm qp}, \Gamma_c}^2(\Omega)$ 
except for a discrete set of wavenumbers $\kappa$. Furthermore, the following estimate holds: 
\begin{eqnarray*}
	\|e\|_{H^2(\Omega)}\lesssim \Theta \|u^i\|_{H^2\left(\Omega^{\rm PML}\right)}.
\end{eqnarray*}
\end{theorem}

\begin{proof}
It is evident from the proof of Theorem \ref{MainResult1} that the variational problem  \eqref{aev} has a unique solution in $H_{{\rm qp}, \Gamma_c}^2(\Omega)$ except for a discrete set of wavenumbers $\kappa$. Furthermore, from the well-posedness of the truncated PML problem \eqref{auvpml} and the proof of Theorem \ref{Step2}, it follows that the truncated PML solution $\tilde{u}$ satisfies the stability estimate
\[
\|\tilde{u}\|_{H^2(\Omega^{\rm PML})}\lesssim \|u^i\|_{H^2\left(\Omega^{\rm PML}\right)}.
\]
By applying the inf-sup condition and trace theorem (cf. \cite[Lemmas 2.2 and 2.3]{BL-IP-2014}), we deduce 
\begin{align*}
\|e\|_{H^2(\Omega)}&\lesssim
\sup\limits_{0\neq v\in H_{qp}^{2}(\Omega^{PML}_k)}
	\frac{|a(e,v)|}{\|v\|_{H^2(\Omega^{\rm PML}_k)}}\\
	&\lesssim \sup\limits_{0\neq v\in H_{qp}^{2}(\Omega^{PML}_k)}
	\frac{1}{\|v\|_{H^2(\Omega^{\rm PML}_k)}}\sum\limits_{k=1}^2
	\left|\int_{\Gamma_k} \left(\overline{v} N_k w_k+\partial_{\nu}\overline{v} M_k w_k\right){\rm d}x\right|\\
&\lesssim \sum\limits_{k=1}^2\left(\|M_k w_k\|_{H^{-1/2}(\Gamma_k)}+ \|N_k w_k\|_{H^{-3/2}(\Gamma_k)}\right)\\
&\lesssim  \Theta \left(\|\tilde{u}-u^i\|_{H^{3/2}(\Gamma_1)}
	+\|\partial_{\nu}\tilde{u}-\partial_{\nu} u^i\|_{H^{1/2}(\Gamma_1)}
	+\|\tilde{u}\|_{H^{3/2}(\Gamma_2)}+\|\partial_{\nu}\tilde{u}\|_{H^{1/2}(\Gamma_2)}\right)\\
&\lesssim \Theta \|u^i\|_{H^2(\Omega^{\rm PML})}.	
\end{align*}
which completes the proof. 
\end{proof}

It can be observed from \eqref{The} and Theorem \ref{Step3} that the PML solution exhibits exponential convergence as the PML parameters $\sigma_1$ and $\sigma_2$, or the PML thickness $\delta$, increase.

\section{Numerical methods}\label{section6}

In this section, drawing inspiration from the mixed finite element method proposed in \cite{YL-JCP-2023} for the cavity scattering problem, we present three different decomposition methods aimed at solving the truncated PML problem.

\subsection{$(q, p)$-decomposition}\label{sub-qp}

Consider the coupled PML problem:
\begin{equation}\label{pqPML}
\left\{
\begin{aligned}
&\tilde{\Delta} p+\kappa^2 p=\frac{1}{2\kappa^2} f\qquad &&{\rm in}\,\Omega^{\rm PML},\\
&\tilde{\Delta} q-\kappa^2 q=\frac{1}{2\kappa^2} f\qquad &&{\rm in}\,\Omega^{\rm PML},\\
&q=p,\quad \partial_{\nu} q=\partial_{\nu} p &&{\rm on}\,\Gamma_2^{\rm PML},\\
&q=p,\quad \partial_{\nu} q=\partial_{\nu} p &&{\rm on}\,\Gamma_c,\\
&q-p=u^i,\quad \partial_{\nu} q-\partial_{\nu} p=\partial_{\nu} u^i &&{\rm on}\,\Gamma_1^{\rm PML},
\end{aligned}
\right.
\end{equation}
where $f$ is defined in \eqref{fpml}.

\begin{lemma}
The boundary value problem \eqref{pqPML} is equivalent to the truncated PML problem \eqref{MBHE}.
\end{lemma}

\begin{proof}
Let $u$ denote the solution of \eqref{MBHE}, and define
\[
p=\frac{1}{2k^2}(\tilde{\Delta}u-\kappa^2 u ),\quad q=\frac{1}{2k^2}(\tilde{\Delta}u+\kappa^2 u).
\]
A straightforward computation yields
\begin{align*}
\tilde{\Delta} p+\kappa^2 p &=\frac{1}{2\kappa^2}
\big[\tilde{\Delta}(\tilde{\Delta} u-\kappa^2 u)+\kappa^2 (\tilde{\Delta} u-\kappa^2 u)\big]\\
&=\frac{1}{2\kappa^2}(\tilde{\Delta}^2 u-\kappa^4 u)=\frac{1}{2\kappa^2}f,\\
\tilde{\Delta} q-\kappa^2 q &=\frac{1}{2\kappa^2}
\big[\tilde{\Delta}(\tilde{\Delta} u+\kappa^2 u)-\kappa^2(\tilde{\Delta} u+\kappa^2 u)\big]\\
&=\frac{1}{2\kappa^2}(\tilde{\Delta}^2 u-\kappa^4 u)=\frac{1}{2\kappa^2}f,
\end{align*}
and
\[
	q-p=u, \quad \partial_{\nu}q-\partial_{\nu}p=\partial_{\nu} u,
\]
indicating that $(p, q)$ constitutes the solution of the boundary value problem \eqref{pqPML}. 

On the other hand, if $(p, q)$ is the solution of \eqref{pqPML}, we denote $v=q-p$, leading to
\[
	\tilde{\Delta} v=\tilde{\Delta}q- \tilde{\Delta}p=\left(\frac{1}{2\kappa^2}f+\kappa^2 q\right)
	-\left(\frac{1}{2\kappa^2}f-\kappa^2 p\right)=\kappa^2(q+p)
\]
and
\begin{align*}
	\tilde{\Delta}^2 v-\kappa^4 v &= \tilde{\Delta}(\kappa^2(q+p))-\kappa^4 (q-p)\\
	&= \kappa^2\big[(\tilde{\Delta}q-\kappa^2 q)+(\tilde{\Delta}p+\kappa^2 p)\big]=f.
\end{align*}
On the boundary of the cavity, the following conditions are satisfied:
\[
	v=q-p=0,\quad \partial_{\nu} v=\partial_{\nu} q-\partial_{\nu}p=0\quad\text{on} \, \Gamma_c, 
\]
Similarly, on the outer boundaries of the PML layers, the boundary conditions are as follows:
\begin{equation*}
\begin{aligned}
&v=q-p=u^i,\quad \partial_{\nu} v=\partial_{\nu} q-\partial_{\nu}p=\partial_{\nu} u^i,\quad 
&&{\rm on}\,\Gamma_1^{\rm PML},\\
&v=q-p=0,\quad \partial_{\nu} v=\partial_{\nu} q-\partial_{\nu}p=0,\quad &&{\rm on}\,\Gamma_2^{\rm PML},
\end{aligned}
\end{equation*}
which confirm that the function $v=q-p$ serves as a solution for \eqref{MBHE}.
\end{proof}

Consider $\mathcal{M}_h$ as a triangulation of $\Omega^{\rm PML}$ such that $\overline{\Omega}^{\rm PML}$ is the union of all triangular elements $K$ in $\mathcal{M}_h$. Let $C_h^I$ and $C_h^B$ represent the collections of interior and boundary edges of the mesh $\mathcal{M}_h$, respectively. Denote by $\mathcal{P}_1$ the finite element space associated with $\mathcal{M}_h$ by using the piecewise linear functions. 

Introduce the discrete space
 \[
S_h=\left\{\phi_h\in C(\overline{\Omega}^{\rm PML}): \phi_h|_{K}\in\mathcal{P}_1(K),\, \forall K\in\mathcal{M}_h\right\},
 \]
and its subspaces $S_h^0$ and $S_h^\Omega$, which exhibit vanishing degrees of freedom on the boundary nodes of  $\Gamma:=\Gamma_c\cup\Gamma_1^{\rm PML}\cup\Gamma_1^{\rm PML}$ and the interior nodes of $\Omega^{\rm PML}$, respectively. It is clear that $S_h=S_h^0\oplus S_h^\Omega$.
 
For any $K, K'\in \mathcal{M}_h$, if $e=\partial K\cap \partial K'$, the jump is defined as
 \[
	\left[\partial_\nu \phi\right]=\partial_{\nu}^{K} u|_{K}-\partial^{K'}_\nu u|_{K'},
 \]
where $\partial_{\nu}^K u|_{K}$ denotes the normal derivative of $u|_K$ on $\partial K$. On the boundary edges $e=\partial K\cup \Gamma$, the jump is set to be zero, i.e., $[\partial_{\nu} u]=0$.

For any functions $\varphi, \psi\in S_h$, the sesquilinear form $G: S_h\times S_h\rightarrow \mathbb{C}$ for the penalty term is defined as follows:
\[
G(\phi, \psi)=\sum\limits_{e\in C^I_h} \eta_e h_e \int_e [\partial_\nu \phi][\partial_\nu \psi]{\rm d}s,
\]
where $h_e$ denotes the length of edge $e$, and the penalty parameter $\eta_e$ is a complex number having a positive real part. 
 
To enforce the boundary conditions, we approximate $p$ and $q$ by 
\[
p\approx w_h+p_h,\quad q\approx w_h+q_h+u^i|_{\Gamma_1^{\rm PML}},\quad p_h, q_h\in S_h^0, \quad w_h\in S_h^{\Omega}.
\]
Multiplying the Helmholtz and modified Helmholtz components by any $\varphi_h\in S_h^0$, we deduce
 \begin{align}
	-\frac{1}{2\kappa^2}\int_{\Omega_1^{\rm PML}} \sigma f\overline{\varphi}_h \,{\rm d}x
	&=b_p(p_h+w_h, \varphi_h)+G(p_h+w_h, \varphi_h),\label{qp1}\\
	-\frac{1}{2\kappa^2}\int_{\Omega_1^{\rm PML}} \sigma f\overline{\varphi}_h  \,{\rm d}x
	&=b_q(q_h+w_h, \varphi_h)-G(q_h+w_h, \varphi_h),\label{qp2}
 \end{align}
 where the sesquilinear form $b_p: S_h\times S_h\rightarrow \mathbb{C}$
 and $b_q: S_h\times S_h\rightarrow \mathbb{C}$  are given by
 \begin{align*}
 b_p(p_h, \varphi_h)&:=\int_{\Omega^{\rm PML}} \left(\frac{1}{\sigma}
 \frac{\partial p_h}{\partial x_2}\frac{\partial \overline{\varphi}_h}{\partial x_2}
	+\sigma\frac{\partial p_h}{\partial x_1}\frac{\partial \overline{\varphi}_h}{\partial x_1}
	-\kappa^2\sigma p_h\overline{\varphi}_h\right){\rm d}x,\\
 b_q(q_h, \varphi_h)&:=\int_{\Omega^{\rm PML}} \left(\frac{1}{\sigma}
 	\frac{\partial q_h}{\partial x_2}\frac{\partial \overline{\varphi}_h}{\partial x_2}
	+\sigma\frac{\partial q_h}{\partial x_1}\frac{\partial \overline{\varphi}_h}{\partial x_1}
	+\kappa^2\sigma q_h\overline{\varphi}_h\right){\rm d}x.
 \end{align*}
 Multiplying the Helmholtz and modified Helmholtz components with $\psi_h\in S_h^\Omega$, and subtracting them from each other yields
 \begin{align}\label{qp3}
-\int_{\Gamma_1^{\rm PML}}\frac{1}{\sigma} \overline{\psi_h} \partial_{x_2 }u^i {\rm d}x_1
 &=b_p(p_h, \psi_h)-b_q(q_h, \psi_h)-2\kappa^2\int_{\Omega^{\rm PML}} \sigma w_h \overline{\psi_h}{\rm d}x\notag\\
&\quad +G(p_h+q_h+2w_h, \psi_h).
 \end{align}
 
 Denote the bases for $S^0_h$ and $S_h^\Omega$ by $\left\{\varphi_i\right\}_{i=1}^{N_I}$ and $\left\{\psi_i\right\}_{i=1}^{N_B}$, respectively. Introduce the related matrices as follows: 
 \begin{align*}
&	M_{pp}(i, j):=b_p(\varphi_j, \varphi_i)+G(\varphi_j, \varphi_i),\quad
	M_{pw}(i, j):=b_p(\psi_j, \varphi_i)+G(\psi_j, \varphi_i),\\
&	M_{qq}(i, j):=b_q(\varphi_j, \varphi_i)+G(\varphi_j, \varphi_i),\quad
	M_{qw}(i, j):=b_q(\psi_j, \varphi_i)+G(\psi_j, \varphi_i),\\
&	M_{wp}(i, j):=b_p(\varphi_j, \psi_i),\quad
	M_{wq}(i, j):=b_q(\varphi_j, \psi_i),\\
&	M_{ww}(i, j)=-2\kappa^2\int_{\Omega^{\rm PML}} \sigma \psi_j \overline{\psi}_i\,{\rm d}x
	+G(\psi_j, \psi_i),
 \end{align*}
and vectors
\begin{eqnarray*}
	F_1(i)=-\frac{1}{2\kappa^2}\int_{\Omega_1^{\rm PML}} \sigma f\overline{\varphi}_i  {\rm d}x,\quad
	F_2(i)=\int_{\Gamma_1^{\rm PML}}\frac{1}{\sigma} \overline{\psi_i}\partial_{x_2 }u^i {\rm d}x_1.
\end{eqnarray*}
The discretized form can be described by
\[
\begin{bmatrix}
M_{pp} & 0 & M_{pw}\\
0 & M_{qq} & M_{qw}\\
M_{wp} & -M_{wq} & M_{ww}
\end{bmatrix}
\begin{bmatrix}
q_h\\ p_h\\ w_h
\end{bmatrix}
=\begin{bmatrix}
F_1\\ F_1 \\ F_2
\end{bmatrix}.
\]
After solving the linear system described above, we can obtain the numerical solution of \eqref{pqPML} as follows: $p=p_h+w_h$ and $q=q_h+w_h+u^i|_{\Gamma_1^{\rm PML}}$. Numerically, we observe that in order to reduce the oscillation of the bending moment, it is important to use the correct penalty terms as given in \eqref{qp1}--\eqref{qp3}.

\subsection{$(u, q)$-decomposition}\label{sub-uq}

For the $(q, p)$-decomposition method as stated in \eqref{pqPML}, the decomposed components $p, q$ are coupled on the boundaries, necessitating a special treatment in the discretization process. If we consider the decomposition $(u, q)$, where $u=q-p$ and $(q, p)$ are the solution of \eqref{pqPML}, then it is straightforward to verify that
\begin{align*}
\tilde{\Delta} u+\kappa^2 u&=\tilde{\Delta}\left(q-p\right)+\kappa^2(q-p)\\
&=\tilde{\Delta}q-\kappa^2 q+2\kappa^2 q-(\tilde{\Delta }p+\kappa^2 p)=2\kappa^2 q.
\end{align*}
Using the boundary conditions of $u$, we obtain 
 \begin{equation}\label{uqPML}
\left\{
\begin{aligned}
&\tilde{\Delta} u+\kappa^2 u-2\kappa^2 q=0\qquad &&{\rm in}\,\Omega^{\rm PML},\\
&\tilde{\Delta} q-\kappa^2 q=\frac{1}{2\kappa^2} f\qquad &&{\rm in}\,\Omega^{\rm PML},\\
&u=0,\quad \partial_\nu u=0 &&{\rm on}\,\Gamma_2^{\rm PML},\\
&u=0,\quad \partial_{\nu} u=0 &&{\rm on}\,\Gamma_c,\\
&u=u^i,\quad \partial_{\nu} u=\partial_{\nu} u^i &&{\rm on}\,\Gamma_1^{\rm PML},
\end{aligned}
\right.
\end{equation}
where $f$ is defined in \eqref{fpml}. It is straightforward to verify that the decomposition \eqref{uqPML} is equivalent to the decomposition \eqref{pqPML}, and both are equivalent to the truncated PML problem \eqref{MBHE}.
 
 Multiplying the modified Helmholtz equation by $\varphi\in H_{{\rm qp}, 00}^1(\Omega)$ and applying integration by parts, we obtain
 \begin{align*}
 - \frac{1}{2\kappa^2}\int_{\Omega_1^{\rm PML}} \sigma f \overline{\varphi}{\rm d}x=
 \int_{\Omega^{\rm PML}} \left(\frac{1}{ \sigma}\frac{\partial q}{\partial x_2}\frac{\overline{\partial \varphi}}{\partial x_2}+ \sigma\frac{\partial q}{\partial x_1}\frac{\overline{\partial \varphi}}{\partial x_1}
	+\kappa^2 \sigma q\overline{\varphi}\right){\rm d}x,
 \end{align*}
 and for $\psi\in H_{\rm qp}^1(\Omega^{\rm PML})$, we have 
  \begin{align*}
 0&= \int_{\Omega^{\rm PML}}  \sigma \big[\tilde{\Delta} u+\kappa^2 u-2\kappa^2 q \big]\psi {\rm d}x\\
 &= -\int_{\Omega^{\rm PML}} \left(\frac{1}{\sigma}
 	\frac{\partial u}{\partial x_2}\frac{\overline{\partial \psi}}{\partial x_2}
	+ \sigma\frac{\partial u}{\partial x_1}\frac{\overline{\partial \psi}}{\partial x_1}
	-\kappa^2 \sigma u\overline{\psi}+2\kappa^2  \sigma q\overline{\psi}\right){\rm d}x\\
&\quad -\int_{\Gamma_c} \partial_\nu u\overline{q}{\rm d}x_1
	+\int_{\Gamma_1} \frac{1}{ \sigma}\partial_\nu u \overline{\psi}{\rm d}x_1
	+\int_{\Gamma_2} \frac{1}{ \sigma}\partial_\nu u \overline{\psi}{\rm d}x_1\\
&=	-\int_{\Omega^{\rm PML}} \left(\frac{1}{ \sigma}
 	\frac{\partial u}{\partial x_2}\frac{\overline{\partial \psi}}{\partial x_2}
	+ \sigma\frac{\partial u}{\partial x_1}\frac{\overline{\partial \psi}}{\partial x_1}
	-\kappa^2 \sigma u\overline{\psi}+2\kappa^2  \sigma q\overline{\psi}\right){\rm d}x
	+\int_{\Gamma_1} \frac{1}{ \sigma}\partial_\nu u^i \overline{\psi}{\rm d}x_1.
 \end{align*}

The variational problem of \eqref{uqPML} is to find  $q\in H_{\rm qp}^1(\Omega^{\rm PML})$ and 
$u\in H_{\rm qp}^{1}(\Omega^{\rm PML})$ satisfying $u=u^i$ on $\Gamma_1^{\rm PML}$ and $u=0$ on $\Gamma_c\cup\Gamma_2^{\rm PML}$ such that
 \begin{align}\label{variuq}
 &\int_{\Omega^{\rm PML}} \left(\frac{1}{ \sigma}
 	\frac{\partial u}{\partial x_2}\frac{\overline{\partial \psi}}{\partial x_2}
	+ \sigma\frac{\partial u}{\partial x_1}\frac{\overline{\partial \psi}}{\partial x_1}
	-\kappa^2 \sigma u\overline{\psi}+2\kappa^2  \sigma q\overline{\psi}\right){\rm d}x
	=\int_{\Gamma_1} \frac{1}{ \sigma}\partial_\nu u^{\rm i} \overline{\psi}{\rm d}x_1,\quad
	\forall\,\psi\in H_{\rm qp}^1(\Omega^{\rm PML}),\notag\\
 &\int_{\Omega^{\rm PML}} \left(\frac{1}{ \sigma}
 	\frac{\partial q}{\partial x_2}\frac{\overline{\partial \varphi}}{\partial x_2}
	+ \sigma\frac{\partial q}{\partial x_1}\frac{\overline{\partial \varphi}}{\partial x_1}
	+\kappa^2 \sigma q\overline{\varphi}\right){\rm d}x
	=- \frac{1}{2\kappa^2}\int_{\Omega_1^{\rm PML}} \sigma f \overline{\varphi}{\rm d}x
	,\quad \forall\,\varphi\in H_{\rm qp, 00}^1(\Omega^{\rm PML}).
 \end{align}

 \begin{lemma}
 If $u, q\in H_{\rm qp}^2(\Omega^{\rm PML})$, then the solution of \eqref{variuq} is also the solution of \eqref{uqPML}.
 \end{lemma}
 
 \begin{proof}
Taking $\varphi\in H_{\rm qp, 00}^{1}(\Omega^{\rm PML})$ in the second equation of \eqref{variuq}
 and using integration by part, we obtain 
 \begin{align*}
 & \int_{\Omega^{\rm PML}} \left(\frac{1}{ \sigma}
 	\frac{\partial q}{\partial x_2}\frac{\overline{\partial \varphi}}{\partial x_2}
	+ \sigma\frac{\partial q}{\partial x_1}\frac{\overline{\partial \varphi}}{\partial x_1}
	+\kappa^2 \sigma q\overline{\varphi}\right){\rm d}x\\
&=\int_{\Omega^{\rm PML}}\left[-
	\frac{1}{ \sigma}\nabla\cdot\left(\left[ \sigma\partial_{x_1} q, \frac{1}{ \sigma}\partial_{x_2}q\right]^\top\right)+\kappa^2 q\right]
	 \sigma\overline{\varphi}\,{\rm d}x	
=-\frac{1}{2\kappa^2}\int_{\Omega^{\rm PML}}  \sigma f\overline{\varphi}\,{\rm d}x,
 \end{align*}
 which implies 
 \[
\tilde{\Delta} q-\kappa^2 q=\frac{1}{2\kappa^2} f\quad \text{in}\,\Omega^{\rm PML}.
 \]
 Taking $\psi\in H_{\rm qp, 00}^{1}(\Omega^{\rm PML})$ in the first equation of \eqref{variuq} leads to
 \begin{align*}
 	0&=\int_{\Omega^{\rm PML}} \left(\frac{1}{ \sigma}
 	\frac{\partial u}{\partial x_2}\frac{\overline{\partial \psi}}{\partial x_2}
	+ \sigma\frac{\partial u}{\partial x_1}\frac{\overline{\partial \psi}}{\partial x_1}
	-\kappa^2 \sigma u\overline{\psi}+2\kappa^2  \sigma q\overline{\psi}\right){\rm d}x
	-\int_{\Gamma_1} \frac{1}{ \sigma}\partial_\nu u^i \overline{\psi}{\rm d}x_1\\
	&=\int_{\Omega^{\rm PML}}\left[-
	\frac{1}{ \sigma}\nabla\cdot\left(\left[ \sigma\partial_{x_1} u, \frac{1}{ \sigma}\partial_{x_2}u\right]^\top\right)-\kappa^2 u
	+2\kappa^2 q\right]\sigma\overline{\psi}{\rm d}x,
 \end{align*}
 which verifies 
  \[
	\tilde{\Delta} u+\kappa^2 q-2\kappa^2 q=0\quad \text{in}\,\Omega^{\rm PML}.
 \]
 
 Taking $\psi\in H_{\rm qp}^1(\Omega)$ such that $\partial_{\nu}\psi=0$ 
 on $\Gamma_c$ or $\Gamma_1^{\rm PML}$ or $\Gamma_2^{\rm PML}$
in the first equation of \eqref{variuq}, we deduce that 
 \begin{align*}
  	0&=\int_{\Omega^{\rm PML}} \left(\frac{1}{ \sigma}
 	\frac{\partial u}{\partial x_2}\frac{\overline{\partial \psi}}{\partial x_2}
	+ \sigma\frac{\partial u}{\partial x_1}\frac{\overline{\partial \psi}}{\partial x_1}
	-\kappa^2 \sigma u\overline{\psi}+2\kappa^2  \sigma q\overline{\psi}\right){\rm d}x
	-\int_{\Gamma_1} \frac{1}{ \sigma}\partial_\nu u^{\rm i} \overline{\psi}{\rm d}x_1\\
	&=\int_{\Omega^{\rm PML}}\left[-
	\frac{1}{ \sigma}\nabla\cdot\left(\left[ \sigma\partial_{x_1} u, \frac{1}{ \sigma}\partial_{x_2}u\right]^\top\right)-\kappa^2 u
	+2\kappa^2 q\right]
	 \sigma\overline{\psi}{\rm d}x\\
&\quad-\int_{\Gamma_c} \partial_\nu u\overline{\psi}\,{\rm d}x_1
	+\int_{\Gamma_1} \frac{1}{ \sigma}\partial_\nu u \overline{\psi}{\rm d}x_1
	+\int_{\Gamma_2} \frac{1}{ \sigma}\partial_\nu u \overline{\psi}{\rm d}x_1
	-\int_{\Gamma_1} \frac{1}{ \sigma}\partial_\nu u^{\rm i} \overline{\psi}{\rm d}x_1,
 \end{align*}
 which leads to the natural boundary conditions
 \[
	\partial_\nu u=\left\{
	\begin{aligned}
	& 0\quad &{\rm on}\,\Gamma_c\cup\Gamma^{\rm PML}_2,\\
	&\partial_\nu u^i &{\rm on}\,\Gamma_1^{\rm PML},
	\end{aligned}
	\right.
 \]
 thus completing the proof. 
 \end{proof}

Taking the same notation as in Subsection \ref{sub-qp}, we can approximate the solution $(u, q)$ as follows:
\[
u\approx u_h,\quad q\approx q_h, \quad u_h\in S_0, \quad q_h\in S_h.
\] 
For any $ \psi_h\in S_h$ and $\varphi_h\in S_h^0$, we solve the following problem: 
 \begin{align*}
 &\int_{\Omega^{\rm PML}} \left(\frac{1}{ \sigma}
 	\frac{\partial u_h}{\partial x_2}\frac{\overline{\partial \psi}_h}{\partial x_2}
	+ \sigma\frac{\partial u_h}{\partial x_1}\frac{\overline{\partial \psi}_h}{\partial x_1}
	-\kappa^2 \sigma u_h\overline{\psi}_h+2\kappa^2  \sigma q_h\overline{\psi}_h\right){\rm d}x\\
&\quad+\sum\limits_{e\in C^I_h} \eta_e h_e 
	\int_e [\partial_\nu u_h]\overline{[\partial_\nu \psi_h]}{\rm d}s
	=\int_{\Gamma_1} \frac{1}{ \sigma}\partial_\nu u^i \overline{\psi}_h{\rm d}x_1,\\
 &\int_{\Omega^{\rm PML}} \left(\frac{1}{ \sigma}
 	\frac{\partial q_h}{\partial x_2}\frac{\overline{\partial \varphi}_h}{\partial x_2}
	+ \sigma\frac{\partial q_h}{\partial x_1}\frac{\overline{\partial \varphi}_h}{\partial x_1}
	+\kappa^2 \sigma q_h\overline{\varphi}_h\right){\rm d}x
	=- \frac{1}{2\kappa^2}\int_{\Omega_1^{\rm PML}} \sigma f \overline{\varphi}_h {\rm d}x. 
 \end{align*}

 \subsection{Decoupled $(p, q)$-decomposition}\label{sub-dqp}
 
It follows from \cite{YL-JCP-2023} that if we take 
\[
	p=\frac{1}{2k^2}\left(\Delta u-\kappa^2 u\right),\quad
	q=\frac{1}{2k^2}\left(\Delta u+\kappa^2 u\right), 
\]
where $u$ is the solution of \eqref{TotalBiharmonicReduced}, then $(p, q)$ satisfies the following boundary value problem:
 \begin{equation}\label{PQSystem}
\left\{
\begin{aligned}
&\Delta p+\kappa^2  p=0,\quad \Delta q-\kappa^2 q=0\qquad &{\rm in}\,\Omega,\\
&p-q=0,\quad \partial_{\nu} p-\partial_{\nu} q=0 &{\rm on}\,\Gamma_c,\\
&\partial_\nu p=T_1 p+g\quad \partial_{\nu} q=-T_2 q &{\rm on}\,\Gamma_1,\\
&\partial_\nu p=T_1 p\quad \partial_{\nu} q=-T_2 q &{\rm on}\,\Gamma_2,
\end{aligned}
\right.
\end{equation}
where $g=4{\rm i}\beta\kappa^2 e^{{\rm i}\left(\alpha x-\beta y\right)}$, and
\[
	T_1 f=\sum\limits_{n\in\mathbb{Z}}{\rm i}\beta_n f_n e^{{\rm i}\,\alpha_n x_1},\quad
	T_2 f=\sum\limits_{n\in\mathbb{Z}}{\rm i}\gamma_n f_n e^{{\rm i}\,\alpha_n x_1},\quad
	f_n=\frac{1}{\Lambda}\int_0^\Lambda fe^{-{\rm i}\alpha_n x_1}{\rm d}x.
\]

It is evident from \eqref{PQSystem} that $p$ and $q$ are decoupled except on $\Gamma_c$, and they satisfy the bounded outgoing wave condition.  This enables us to approximate \eqref{PQSystem} using the PML technique:
\begin{equation*}
\left\{
\begin{aligned}
	&\frac{\partial }{\partial x_1}\left(\sigma\frac{\partial p}{\partial x_1}\right)
	+\frac{\partial }{\partial x_2}\left(\frac{1}{\sigma}\frac{\partial p}{\partial x_2}\right)
	+\kappa^2 \sigma p=-f\quad&&{\rm in}\,\Omega^{\rm PML},\\
	&\frac{\partial }{\partial x_1}\left(\sigma\frac{\partial q}{\partial x_1}\right)
	+\frac{\partial }{\partial x_2}\left(\frac{1}{\sigma}\frac{\partial q}{\partial x_2}\right)
	-\kappa^2 \sigma q=0\quad&&{\rm in}\,\Omega^{\rm PML},\\
	&p=q,\quad \partial_\nu p=\partial_\nu q &&{\rm on}\,\Gamma_c,\\
	&p=-2\kappa^2 u^i,\quad q=0 &&{\rm on}\,\Gamma^{\rm PML}_1,\\
	&p=0,\quad q=0 &&{\rm on}\,\Gamma^{\rm PML}_2
\end{aligned}
\right.
\end{equation*}
where
\begin{eqnarray*}
f(x)=
\left\{
	\begin{aligned}
	& 2\kappa^2\left[\frac{\partial }{\partial x_1}\left(\sigma\frac{\partial u^i}{\partial x_1}\right)+\frac{\partial }{\partial x_2}\left(\frac{1}{\sigma}\frac{\partial u^i}{\partial x_2}\right)
	+\kappa^2 \sigma u^i\right], \quad && x\in \Omega_1,\\
	&0,\quad  &&  x\in \Omega^{\rm PML}\setminus\overline{\Omega_1}.
	\end{aligned}
	\right.
\end{eqnarray*}

We approximate the solution $(p, q)$ as follows: 
 \[
	p\approx w_h+p_h,\quad q\approx w_h+q_h,\quad p_h, q_h\in \hat{S}_h, \quad w_h\in S_h^{\Gamma_c},
 \]
 where the spaces $\hat{S}_h$ and $S_h^{\Gamma_c} $ consist of functions with vanishing degrees of freedom on the nodes of $\Gamma_c$ and elsewhere, respectively. For any $\varphi_h\in \hat{S}_h$ and $\psi_h\in S_h^{\Gamma_c}$, the discrete problem is to find $p_h=-2\kappa^2 u^i$ on $\Gamma_c$ such that
 \[
 \left\{
 \begin{aligned}
	&b_p(p_h+w_h, \varphi_h)+G(p_h+w_h, \varphi_h)=\int_{\Omega_1^{\rm PML}} \sigma f\overline{\varphi}_h \,{\rm d}x,\\
	&b_q(q_h+w_h, \varphi_h)-G(q_h+w_h, \varphi_h)=0,\\
	&b_p(p_h, \psi_h)-b_q(q_h, \psi_h)
	-2\kappa^2\int_{\Omega^{\rm PML}} \sigma w_h \overline{\psi_h}\,{\rm d}x+G(p_h+q_h+2w_h, \psi_h)=0.
 \end{aligned}
 \right.
 \]
 
Compared to the $(q, p)$-decomposition, this decoupled $(q, p)$-decomposition can handle the boundary condition on $\Gamma_1$ much more easily. In the numerical experiments, we observe that the result obtained from the decoupled $(q, p)$-decomposition is more stable than the one obtained from the $(q, p)$-decomposition.

\subsection{Numerical experiments}

In this subsection, we present some numerical experiments to demonstrate the effectiveness of the mixed finite element methods outlined in Subsections \ref{sub-qp}--\ref{sub-dqp}.

In the following examples, we consider an incident plane wave with wavenumber $\kappa=\pi$, an incident angle $\theta=\pi/3$, and set the PML layers beyond $\Gamma_1=\left\{x_2=0.5\right\}$ and $\Gamma_2=\left\{x_2=-0.5\right\}$ with a thickness of $\Delta h=2.5$. We choose parameters $m=4$, $\sigma_1=14$, and $\sigma_2=5$ to satisfy the condition \eqref{PMLrestriction}. The penalty parameter is set to a complex constant value $\eta=0.001+0.001{\rm i}$. The numerical convergence rate is obtained by comparing the results obtained using a uniform mesh $h=[0.015:0.005:0.05]$ with the reference field obtained using a refined mesh $h=0.008$ based on the decoupled $(q, p)$-method.

\subsubsection{Example 1}

\begin{figure}
\centering
\includegraphics[width=0.45\textwidth]{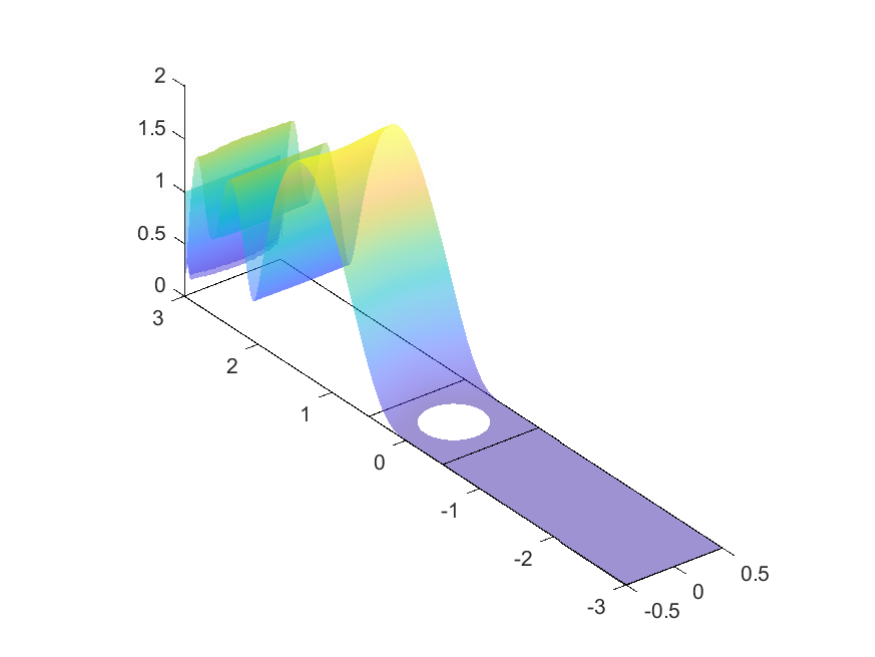}
\includegraphics[width=0.45\textwidth]{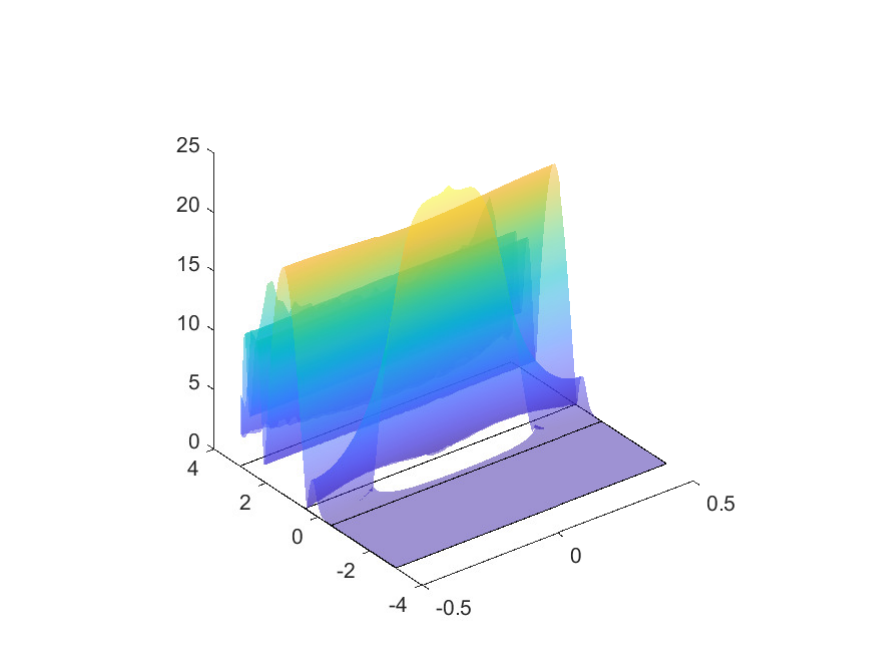}
\caption{The solution of $(q, p)$-decomposition for a circle-shaped cavity.}
\label{qp-circle}
\end{figure}

\begin{figure}
\centering
\includegraphics[width=0.45\textwidth]{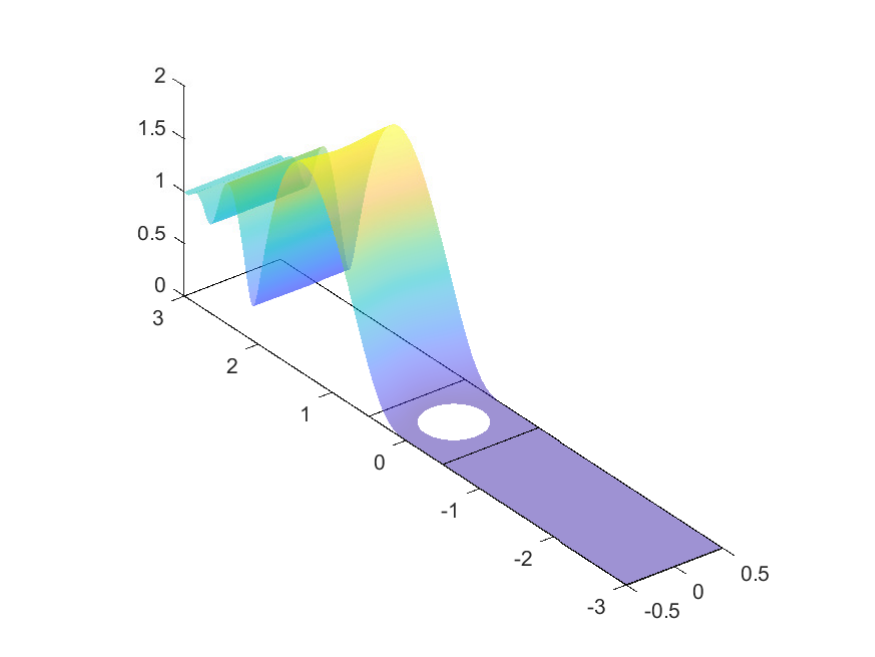}
\includegraphics[width=0.45\textwidth]{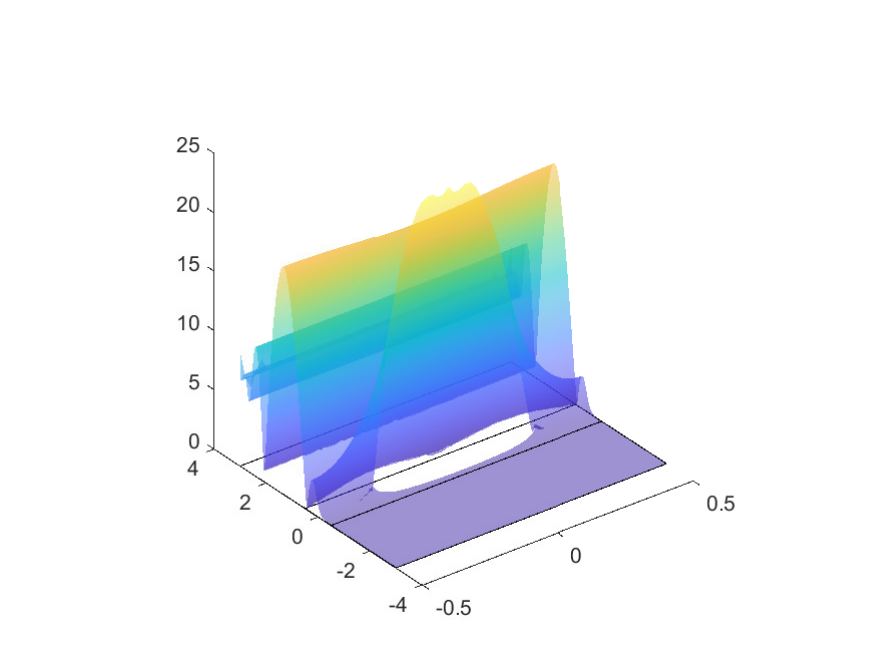}
\caption{The solution of $(u, q)$-decomposition for a circle-shaped cavity.}
\label{uq-circle}
\end{figure}

\begin{figure}
\centering
\includegraphics[width=0.45\textwidth]{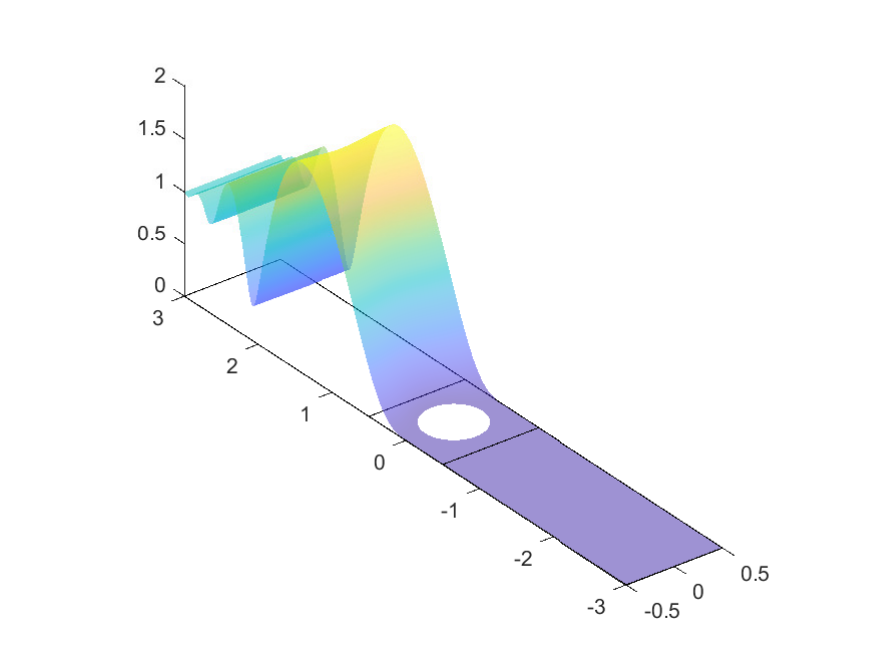}
\includegraphics[width=0.45\textwidth]{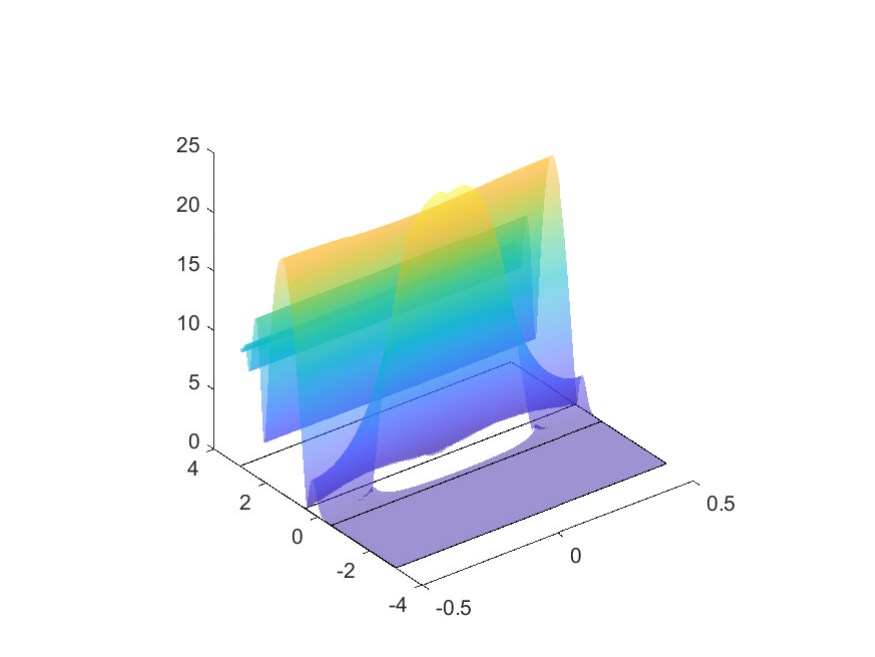}
\caption{The solution of decoupled $(q, p)$-decomposition for a circle-shaped cavity.}
\label{dqp-circle}
\end{figure}

In the first example, the cavity is a disk with a radius $R=0.3$. Figures \ref{qp-circle}--\ref{dqp-circle} display the results obtained from the three methods, illustrating the influence of the interior penalty. The mesh size is chosen as $h=0.02$ with a number of nodal points 16090 in the computational domain. In each figure, the left and right panels represent the $L^2$-norm of the field $u$ and its bending moment $\Delta u$, respectively. It is evident that the field decays exponentially in $\Omega_2$ using all three methods. When compared with the results in \cite{YL-JCP-2023}, there is no oscillation near the surface of the cavity. This illustrate the necessity of applying the interior penalty to suppress the oscillation of the bending moment in practice.

\begin{figure}
\centering
\includegraphics[width=0.45\textwidth]{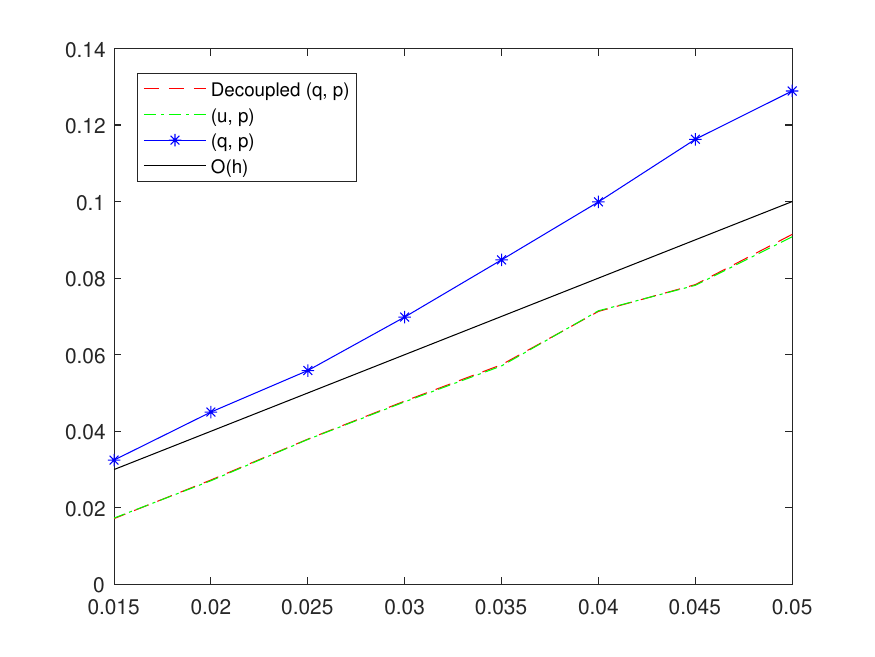}
\includegraphics[width=0.45\textwidth]{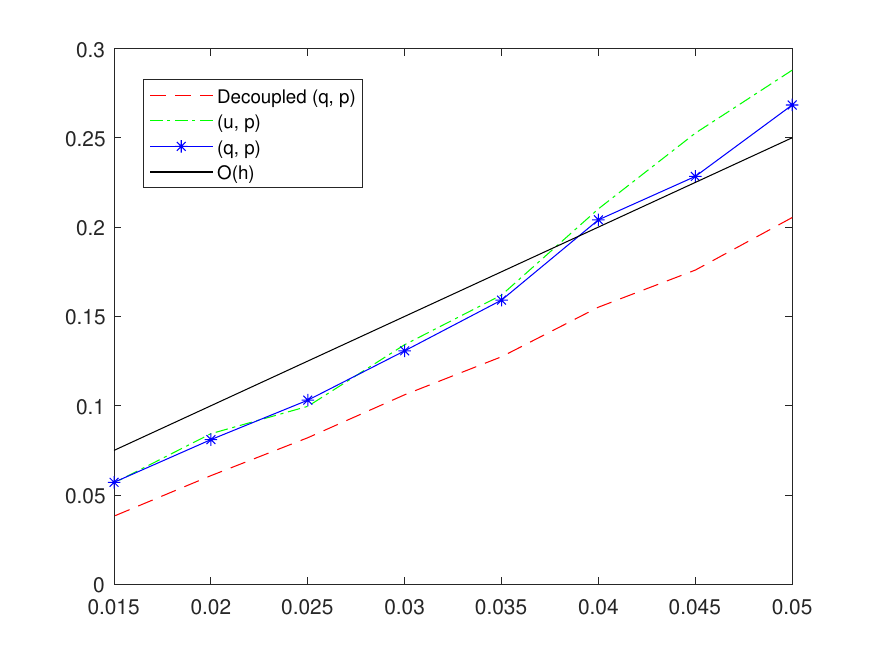}
\caption{The convergent rate of the three methods for a circle-shaped cavity.}
\label{rate-circle}
\end{figure}

We also assess the convergence rate with respect to the mesh size. Figure \ref{rate-circle} illustrates the convergence rate for $u$ and $\Delta u$, respectively. In both figures, the red-dashed line, green-dashed-point line, and blue-star line represent the results for the decoupled $(q, p)$-, $(u, q)$-, and $(q, p)$-decomposition, respectively. The black-solid line serves as the reference curve indicating first order convergence. When the mesh size is sufficiently small, the convergence appears to be approximately linear. It is important to note that when the mesh is coarse, the result obtained from $(q, p)$-decomposition is unsatisfactory. This may be attributed to the inherent coupled boundary condition on $\Gamma^{\rm PML}_1$.

\subsubsection{Example 2} 
\begin{figure}
\centering
\includegraphics[width=0.45\textwidth]{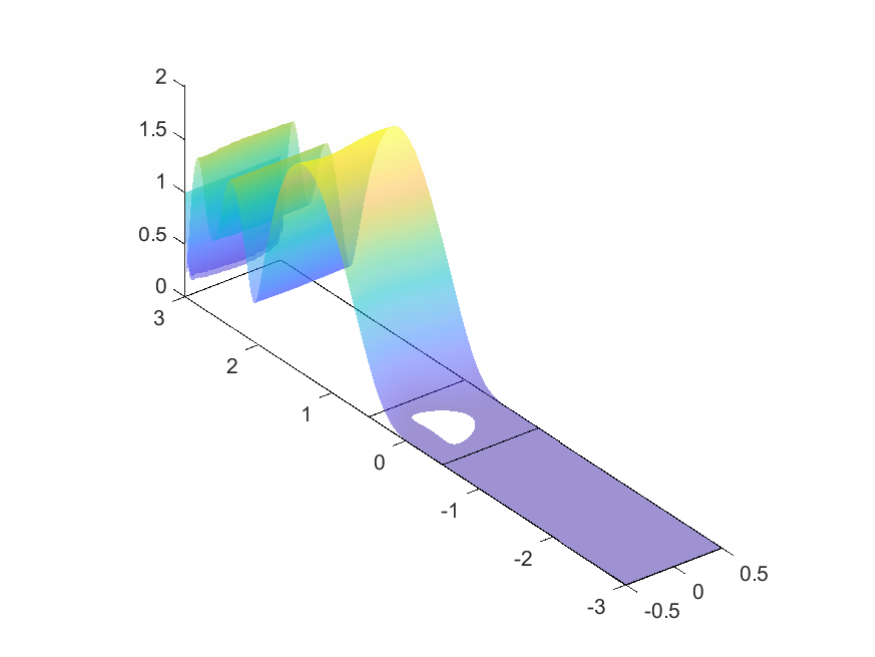}
\includegraphics[width=0.45\textwidth]{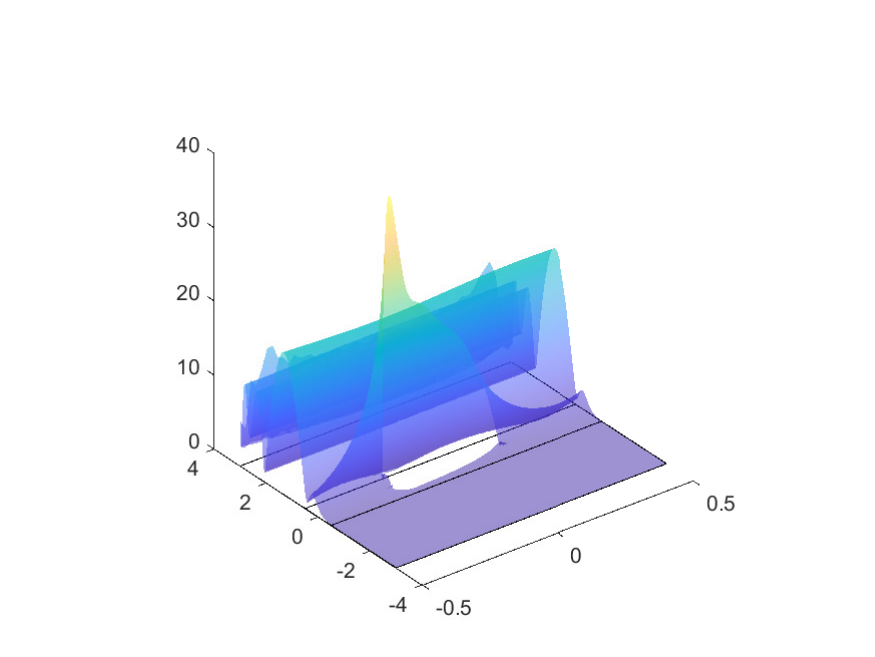}
\caption{The solution of $(q, p)$-decomposition for a kite-shaped cavity.}
\label{qp-kite}
\end{figure}

\begin{figure}
\centering
\includegraphics[width=0.45\textwidth]{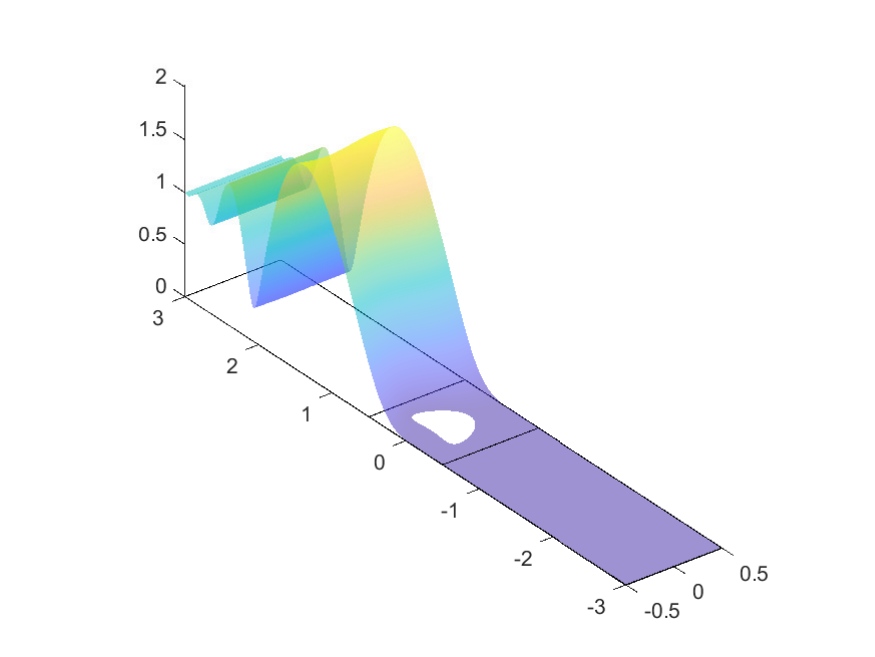}
\includegraphics[width=0.45\textwidth]{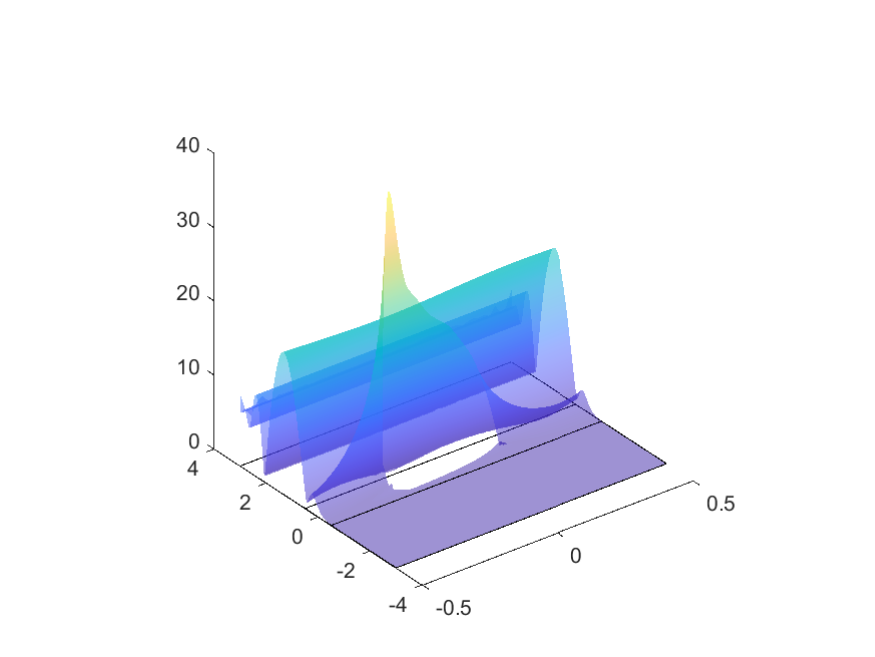}
\caption{The solution of $(u, q)$-decomposition for a kite-shaped cavity.}
\label{uq-kite}
\end{figure}

\begin{figure}
\centering
\includegraphics[width=0.45\textwidth]{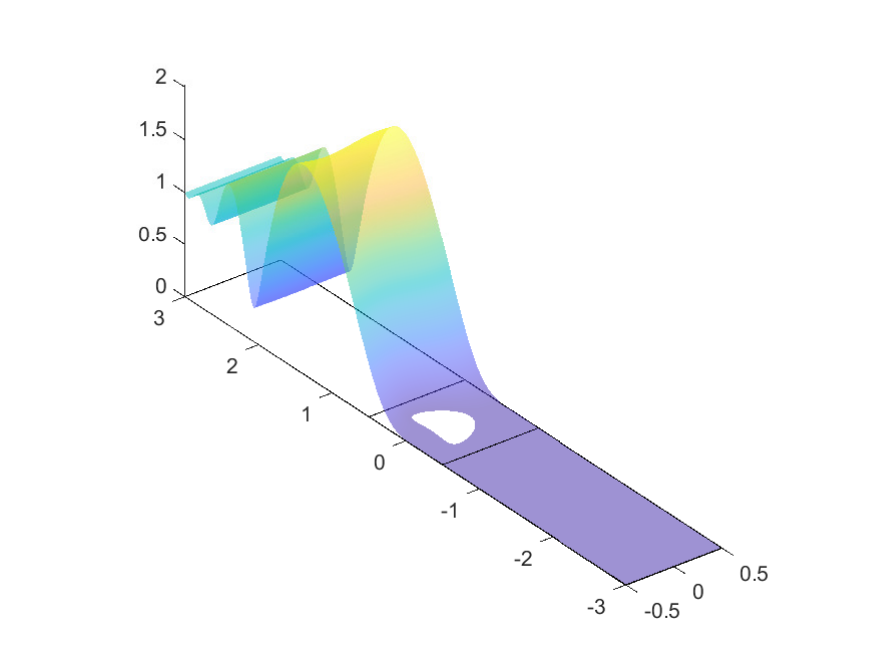}
\includegraphics[width=0.45\textwidth]{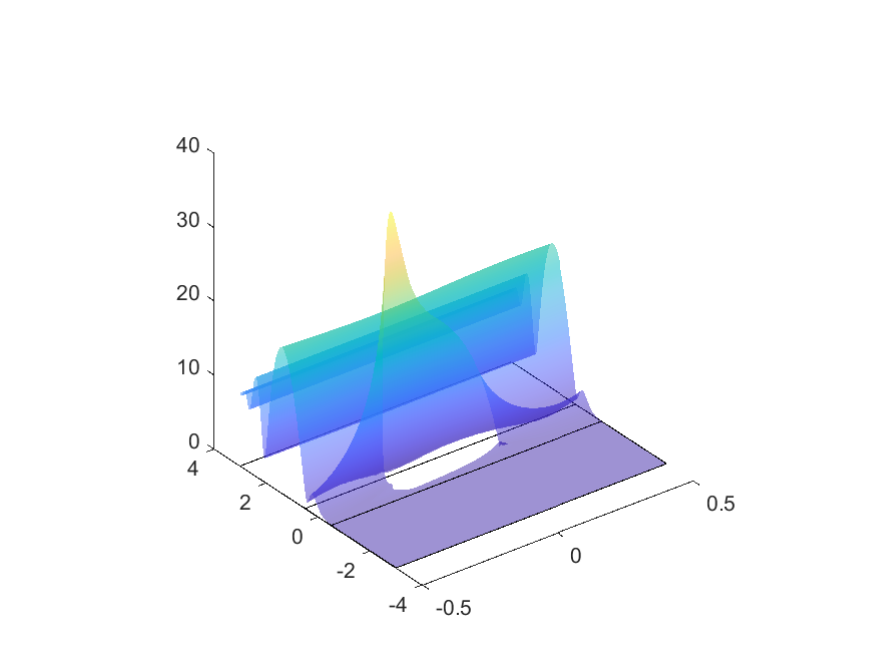}
\caption{The solution of decoupled $(q, p)$-decomposition for a kite-shaped cavity.}
\label{dqp-kite}
\end{figure}

In the second example, we examine a kite-shaped cavity defined by 
\[
	x(t)=0.2\left(\cos t+0.07\cos 2t-0.1,  0.3\sin t\right),\quad 0\leq t\leq 2\pi.
\]
The other parameters remain consistent with Example 1. Figures \ref{qp-kite}--\ref{dqp-kite} depict the results and illustrate the impact of the interior penalty for the three methods. The mesh size is set to $h=0.02$ with a total of 17,804 nodal points in the computational domain. The $L^2$-norm of the field $u$ and its bending moment $\Delta u$ are presented in the left and right panels, respectively. Similarly, it is evident that the field exponentially decays in $\Omega_2$ for all three methods, confirming the effectiveness of the PML method. The addition of the interior penalty term effectively suppresses oscillations near $\Gamma_c$.

\begin{figure}
\centering
\includegraphics[width=0.45\textwidth]{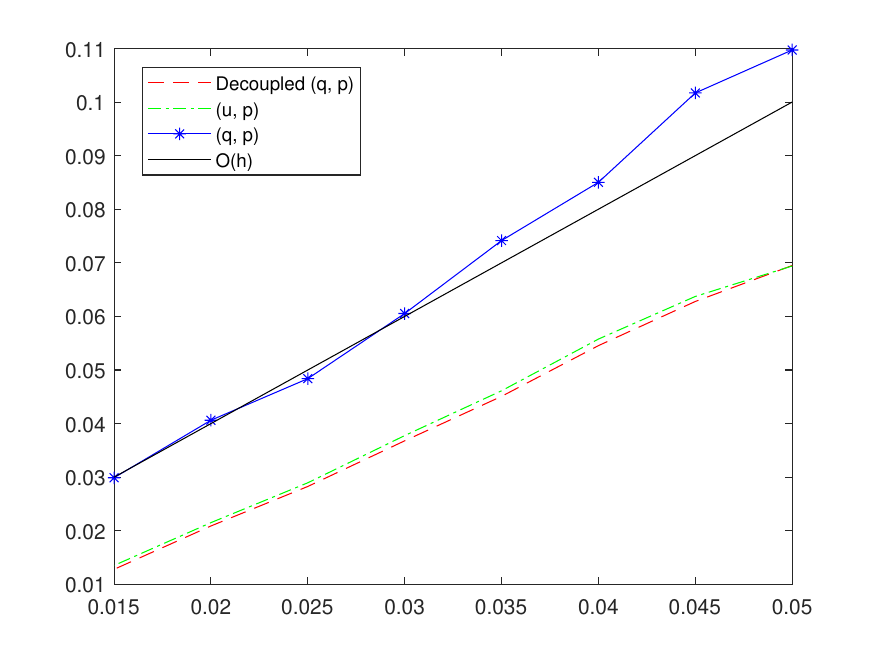}
\includegraphics[width=0.45\textwidth]{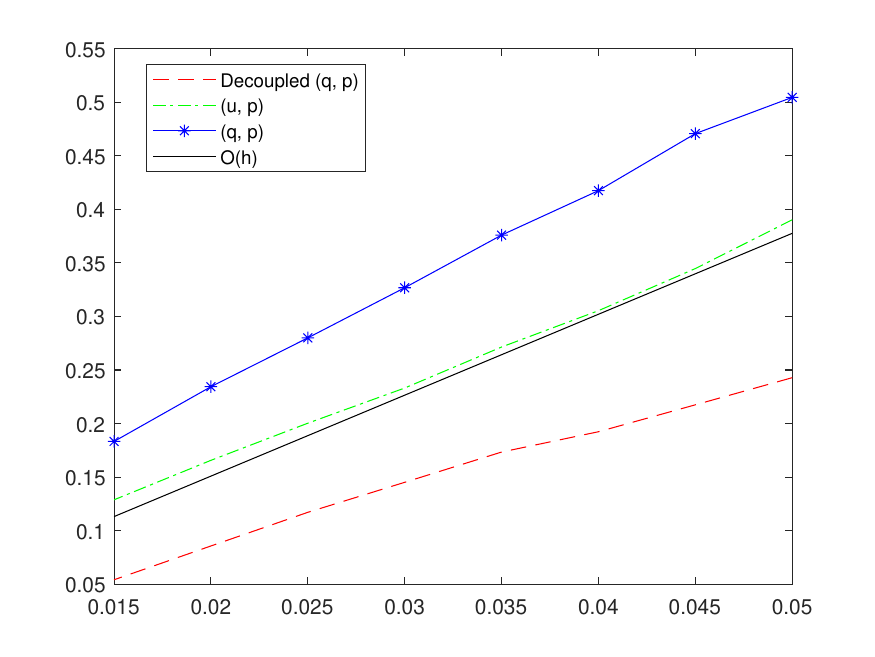}
\caption{The convergent rate of the three methods for a kite-shaped cavity.}
\label{rate-kite}
\end{figure}

Figure \ref{rate-kite} depicts the convergence order for $u$ and $\Delta u$, respectively. In both figures, the red-dashed line, green-dashed-point line, and blue-star line represent the results for the decoupled $(q, p)$-, $(u, q)$-, and $(q, p)$-decomposition, respectively. Again, the black-solid line serves as the reference curve indicating first order convergence. When the mesh size is sufficiently small, the convergence appears to be approximately linear. Additionally, it is observed that the $(q, p)$-decomposition does not perform as well as the other two methods with a coarse mesh.  

\section{Conclusion}\label{section7}

In this paper, we have examined the PML problem concerning the scattering of biharmonic waves in a periodic structure. The well-posedness of the associated variational problem is established in the computational domain. The exponential convergence is achieved for the PML solution with respect to the PML parameter and the thickness of the PML layer. Three mixed finite element methods with interior penalty technique are proposed for solving the PML problem. Numerical experiments are presented to demonstrate the efficiency of the proposed methods. 

In comparison with the method outlined in \cite{BLY-2023-sub}, the current study does not necessitate the explicit derivation of the intricate PML-DtN operator. Consequently, our approach holds applicability to many other complex scattering model problems. Specifically, we are exploring the application of the PML method to the scattering of biharmonic waves by bounded cavities, a problem previously examined in \cite{DL-2023-arXiv} utilizing the boundary integral equation method. We intend to report our progress in a future publication.

\section*{Appendix: Auxiliary functions}\label{appendix}

Denote the domain $\Omega:=[0, \Lambda]\times[0, h]$, and its boundaries $\Gamma=\left\{x_2=h\right\}$ and $S=\{x_2=0\}$. The purpose of the appendix is to explicitly construct two auxiliary functions that satisfy the following conditions:
\begin{itemize}

\item[(1)] For any given $\varphi\in H^{3/2}(\Gamma)$, construct a function $u\in H_{\rm qp}^2(\Omega)$ satisfying the boundary conditions $u(x_1, h)=\varphi(x_1), u(x_1, 0)=0, \partial_{x_2} u(x_1, h)=0, \partial_{x_2} u(x_1, 0)=0$, and the stability estimate
\[
\|u\|_{H^2(\Omega)}\leq C\|\varphi\|_{H^{3/2}(\Gamma)}.
\]

\item[(2)] For any given $\varphi\in H^{1/2}(\Gamma)$, construct $u\in H_{\rm qp}^2(\Omega)$ satisfying the boundary conditions $u(x_1, h)=0, u(x_1, 0)=0, \partial_{x_2} u(x_1, h)=\varphi(x_1), \partial_{x_2} u(x_1, 0)=0$, and the stability estimate
\[
\|u\|_{H^2(\Omega)}\leq C\|\varphi\|_{H^{1/2}(\Gamma)}.
\]
\end{itemize}

\subsection{Case 1}

 Let $\varphi\in H_{\rm qp}^{3/2}(\Gamma)$ admit the Fourier series expansion 
\[
\varphi(x_1)=\sum\limits_{n\in\mathbb{Z}} \varphi_n e^{{\rm i}\alpha_n x_1},\quad \|\varphi\|_{H^{3/2}(\Gamma)}^2=\sum\limits_{n\in\mathbb{Z}} (1+n^2)^{3/2}|\varphi_n|^2.
\]
Define a function 
\begin{equation}\label{uCase1}
	u(x_1, x_2)=L_1(x_2)\sum\limits_{n\in\mathbb{Z}} g_n(x_2)\varphi_n e^{{\rm i}\alpha_n x_1},
\end{equation}
where $L_1$ is the Hermite interpolation basis function, given by 
\[
L_1(x_2)=\left(3h-2x_2\right)x_2^2\frac{1}{h^3},\quad L_1'(x_2)=\frac{6x_2}{h^2}-\frac{6x_2^2}{h^3},
\]
and 
\[
	g_n(x_2)=e^{-|n|(h-x_2)}+\psi_n(x_2),\quad
	\psi_n(x_2)=\left\{
	\begin{aligned}
	&\sin |n|(h-x_2), & h_n <x_2 \leq h,\\
	&0, &\text{otherwise}. 
	\end{aligned}
	\right.
\]
It is easy to check that
$
	g_n(h)=1,
	g'_n(h)=0, |g_n(x_2)|\leq 2, 
$ and
\begin{eqnarray}
&g'_n(x_2)=|n|e^{-|n|(h-x_2)}+\psi'_n(x_2),\quad
	\psi'_n(x_2)=\left\{
	\begin{aligned}
	&-|n|\cos |n|(h-x_2), & h_n <x_2 \leq h,\\
	&0, &\text{otherwise},
	\end{aligned}
	\right. \label{dgn}\\
&g''_n(x_2)=|n|^2e^{-|n|(h-x_2)}+\psi''_n(x_2),\quad
	\psi''_n(x_2)=\left\{
	\begin{aligned}
	&-|n|^2\sin |n|(h-x_2), & h_n <x_2 \leq h,\\
	&0, &\text{otherwise},
	\end{aligned}
	\right.	\label{ddgn}
\end{eqnarray}
where $h_n=\max(h-\frac{2\pi}{|n|}, 0) $.

Taking the derivative of \eqref{uCase1}, we get 
\begin{eqnarray*}
\partial_{x_2} u=\left(\frac{6x_2}{h^2}-\frac{6x_2^2}{h^3}\right)\sum\limits_{n\in\mathbb{Z}} g_n(x_2)\varphi_n e^{{\rm i}\alpha_n x_1}+\left(\frac{3x_2^2}{h^2}-\frac{2x_2^3}{h^3}\right)\sum\limits_{n\in\mathbb{Z}}g'_n(x_2)\varphi_n e^{{\rm i}\alpha_n x_1},
\end{eqnarray*}
which implies 
\[
\left\{
\begin{aligned}
	&u(x_1, h)=\varphi(x_1),\\
	&u(x_1, 0)=0,
\end{aligned}
\right.
\quad
\left\{
\begin{aligned}	 
\partial_y u(x_1, h)=0,\\
\partial_y u(x_1, 0)=0. 
\end{aligned}
\right.
\]

A straightforward computation yields 
\begin{align*}
\|u\|_{L^2(\Omega)}^2 &= \int_0^{\Lambda}\int_0^h
	 \left|L_1(x_2)\sum\limits_{n\in\mathbb{Z}} e^{-|n|(h-x_2)^2}\varphi_n e^{{\rm i}\alpha_n x_1}
	\right|^2{\rm d}x\\
&\leq  \int_0^{\Lambda}\int_0^h\left|\sum\limits_{n\in\mathbb{Z}} \varphi_n e^{{\rm i}\alpha_n x_1}\right|^2{\rm d}x
\leq C\sum\limits_{n\in\mathbb{Z}} |\varphi_n|^2
\leq C\|\varphi\|_{L^2(\Gamma)}^2,
\end{align*}
\begin{align*}
\|\partial_{x_1} u\|_{L^2(\Omega)}^2 &= \int_0^{\Lambda}\int_0^h
	 \left|L_1(x_2)\sum\limits_{n\in\mathbb{Z}} {\rm i}\alpha_n g_n(x_2)\varphi_n e^{{\rm i}\alpha_n x_1}
	\right|^2{\rm d}x\\
	&\leq C\int_0^{\Lambda}\int_0^h
	 \left|\sum\limits_{n\in\mathbb{Z}} n \varphi_n e^{{\rm i}\alpha_n x_1}\right|^2{\rm d}x
	 \leq C\sum\limits_{n\in\mathbb{Z}} \left(1+n^2\right) |\varphi_n|^2
	 \leq C\|\varphi\|_{H^1(\Gamma)}^2,
\end{align*}
and
\begin{align*}
\|\partial_{x_2} u\|_{L^2(\Omega)}^2 &= \int_0^{\Lambda}\int_0^h \left| L_1'(x_2)\sum\limits_{n\in\mathbb{Z}}g_n(x_2)\varphi_n e^{{\rm i}\alpha_n x_1} +L_1(x_2)\sum\limits_{n\in\mathbb{Z}} g'_n(x_2)\varphi_n e^{{\rm i}\alpha_n x_1}
\right|^2 {\rm d}x\\
&\leq C\int_0^{\Lambda}\int_0^h \left|\sum\limits_{n\in\mathbb{Z}} \varphi_n e^{{\rm i}\alpha_n x_1}
\right|^2 + \left| \sum\limits_{n\in\mathbb{Z}} |n|\varphi_n e^{{\rm i}\alpha_n x_1}
\right|^2 {\rm d}x\\
&\leq C\sum\limits_{n\in\mathbb{Z}}\left( |\varphi_n|^2+n^2 |\varphi_n|^2\right)
\leq C\|u\|_{H^1(\Gamma)}^2.
\end{align*}

Taking the second derivatives of $u$ gives 
\begin{align*}
\left\|\frac{\partial^2 u}{\partial x_1^2}\right\|_{L^2(\Omega)}^2 
&=\lim\limits_{\epsilon\rightarrow 0} \int_0^{\Lambda}\int_0^{h-\epsilon}
	 \left|L_1(x_2)\sum\limits_{n\in\mathbb{Z}} -\alpha^2_n g_n(x_2)\varphi_n e^{{\rm i}\alpha_n x_1}
	\right|^2{\rm d}x\\
&\leq \lim\limits_{\epsilon\rightarrow 0} \int_0^{\Lambda}\int_0^{h-\epsilon}
	 \left|L_1(x_2)\sum\limits_{n\in\mathbb{Z}} -\alpha^2_n e^{-2|n|(h-x_2)^2}\varphi_n e^{{\rm i}\alpha_n x_1}
	\right|^2{\rm d}x\\
&\quad +\lim\limits_{\epsilon\rightarrow 0} \int_0^{\Lambda}\int_{0}^{h-\epsilon}
	 \left|L_1(x_2)\sum\limits_{n\in\mathbb{Z}} -\alpha^2_n \psi_n(x_2)\varphi_n e^{{\rm i}\alpha_n x_1}
	\right|^2{\rm d}x\\
&\leq C\lim\limits_{\epsilon\rightarrow 0}
\sum\limits_{n\in\mathbb{Z}} |\alpha^2_n|^2 |\varphi_n|^2 \left(\int_0^{h-\epsilon}
	 e^{-2|n|(h-x_2)}{\rm d}x_2+\int_{h_n}^{h-\epsilon} \sin^2(|n|(h-x_2)){\rm d}x_2\right)\\
&\leq C \sum\limits_{n\in\mathbb{Z}} \frac{|\alpha^2_n|^2}{|n|} |\varphi_n|^2\leq C\|\varphi\|_{H^{3/2}(\Gamma)}^2.
\end{align*}
Using \eqref{dgn} and \eqref{ddgn} leads to 
\begin{align*}
\left\|\frac{\partial^2 u}{\partial x_1\partial x_2}\right\|_{L^2(\Omega)}^2 &=\lim\limits_{\epsilon\rightarrow 0} \int_0^{\Lambda}\int_0^{h-\epsilon}\left|\sum\limits_{n\in\mathbb{Z}} {\rm i}\alpha_n \left(L'_1(x_2)g_n(x_2)+L_1(x_2)g'_n(x_2)\right)\varphi_n e^{{\rm i}\alpha_n x_1}\right|^2{\rm d}x\\
&\leq C	\lim\limits_{\epsilon\rightarrow 0} \sum\limits_{n\in\mathbb{Z}}|\alpha_n|^2|\varphi_n|^2 \bigg(
\int_0^{h-\epsilon} |n|^2e^{-2|n|(h-x_2)}{\rm d}x_2\\
&\quad +\int_{h_n}^{h-\epsilon} |n|^2\cos^2(|n|(h-x_2)){\rm d}x_2\bigg)\\
&\leq C \lim\limits_{\epsilon\rightarrow 0} \sum\limits_{n\in\mathbb{Z}}|n|^3|\varphi_n|^2=C\|\varphi\|_{H^{3/2}(\Gamma)}^2
\end{align*}
and
\begin{align*}
&\left\|\frac{\partial^2 u}{\partial x_2^2}\right\|_{L^2(\Omega)}^2  =\lim\limits_{\epsilon\rightarrow 0} \int_0^{\Lambda}\int_0^{h-\epsilon}\left|\sum\limits_{n\in\mathbb{Z}} \left(L''_1(x_2)g_n(x_2)+2L'_1(x_2)g'_n(x_2)+g'_n(x_2)\right)\varphi_n e^{{\rm i}\alpha_n x_1}\right|^2{\rm d}x\\
&\leq C\lim\limits_{\epsilon\rightarrow 0} \int_0^{\Lambda}\int_0^{h-\epsilon}\left|\sum\limits_{n\in\mathbb{Z}} g_n(x_2)\varphi_n e^{{\rm i}\alpha_n x_1}\right|^2+ \left|\sum\limits_{n\in\mathbb{Z}}g'_n(x_2)\varphi_n e^{{\rm i}\alpha_n x_1} \right|^2 + \left|\sum\limits_{n\in\mathbb{Z}}g''_n(x_2)\varphi_n e^{{\rm i}\alpha_n x_1}\right|^2 
{\rm d}x\\
&\lesssim \|\varphi\|_{L^2(\Gamma)}^2+\|\varphi\|_{H^1(\Gamma)}^2 +\lim\limits_{\epsilon\rightarrow 0} 
\sum\limits_{n\in\mathbb{Z}}|\varphi_n|^2 \bigg(\int_0^{h-\epsilon} |n|^4e^{-2|n|(h-x_2)}{\rm d}x_2\\
&\hspace{6cm} +\int_{h_n}^{h-\epsilon} |n|^4\sin^2(|n|(h-x_2)){\rm d}x_2\bigg)\\
&\lesssim \|\varphi\|_{L^2(\Gamma)}^2+\|\varphi\|_{H^1(\Gamma)}^2+\|\varphi\|_{H^{3/2}(\Gamma)}^2
\lesssim\|\varphi\|_{H^{3/2}(\Gamma)}^2.
\end{align*}
Combining the above estimates, we show that the function $u$ defined in \eqref{uCase1} satisfies the stability estimate
\[
	\|u\|_{H^2(\Omega)}\leq C\|\varphi\|_{H^{3/2}(\Gamma)}.
\]

Following the same argument, we construct a function
\[
u(x_1, x_2)=L_2(x_2)\sum\limits_{n\in\mathbb{Z}} g_n(x_2)\varphi_n e^{{\rm i}\alpha_n x_1},
\]
where $L_2(x_2)$ is the Hermite interpolation basis function, given by  
\[
L_2(x_2)=\left(h+2x_2\right)\frac{(x_2-h)^2}{h^3},
\]
and
\[
	g_n(x_2)=e^{-|n| x_2}+\psi_n(x_2),\quad
	\psi_n(x_2)=\left\{
	\begin{aligned}
	&\sin |n| x_2, & 0 <x_2\leq \min\left(h, \frac{2\pi}{|n|}\right),\\
	&0, &\text{otherwise}.
	\end{aligned}
	\right.
\]
It is easy to verify that 
\[
\left\{
\begin{aligned}
	&u(x_1, h)=0,\\
	&u(x_1, 0)=\varphi(x_1),
\end{aligned}
\right.
\qquad
\left\{
\begin{aligned}	 
\partial_{x_2} u(x_1, h)=0,\\
\partial_{x_2} u(x_1, 0)=0,
\end{aligned}
\right.
\]
and
\[
\|u\|_{H^2(\Omega)}\leq C\|\varphi\|_{H^{3/2}(S)}.
\]

\subsection{Case 2}

Let $\varphi\in H_{\rm qp}^{1/2}(\Gamma)$ admit the Fourier series expansion 
\[
\varphi(x_1)=\sum\limits_{n\in\mathbb{Z}} \varphi_n e^{{\rm i}\alpha_n x_1},\quad \|\varphi\|_{H^{1/2}(\Gamma)}^2=\sum\limits_{n\in\mathbb{Z}} (1+n^2)^{1/2}|\varphi_n|^2.
\]

Consider an auxiliary problem
\begin{equation}\label{auxiliary}
\left\{
\begin{aligned}
&\Delta^2 u=0\qquad&{\rm in}\,\Omega,\\
&u=0,\quad\partial_{\nu} u=0 &{\rm on}\,S,\\
& u=0,\quad\partial_{\nu }u=\varphi &{\rm on}\,\Gamma.
\end{aligned}
\right.
\end{equation}

\begin{lemma}\label{lemmaau}
The auxiliary problem \eqref{auxiliary} admits a unique weak solution.
\end{lemma}

\begin{proof}
Similar to the previous discussion, for any given $\varphi\in H^{1/2}(\Gamma)$, the variational problem for \eqref{auxiliary} is to find $u\in H_{\rm qp}^{2}(\Omega)$ satisfying the boundary conditions $u=0, \partial_{\nu} u=0$ on $S$ and $u=0, \partial_{\nu} u=\varphi$ on $\Gamma$ such that
\begin{equation}\label{Au}
	A(u, v)=0,\quad\forall v\in H_0^2(\Omega),
\end{equation}
where the sesquilinear form $A(u, v):  H_{{\rm qp}}^2(\Omega)\times  
 H_{{\rm qp}}^2(\Omega)\rightarrow \mathbb{C}$ is given by 
\[
	A(u, v)=
	\int_{\Omega} \bigg[\mu\Delta u\Delta \bar{v}+(1-\mu)\sum\limits_{i,j=1}^2
	 \frac{\partial^2 u}{\partial x_i\partial x_j}
\frac{\partial^2 \bar{v}}{\partial x_i\partial x_j}\bigg]{\rm d}x.
\]
It can be shown that if $\mu\in(0, 1)$, the sesquilinear form satisfies G\r{a}rding's inequality. 

To demonstrate that the problem \eqref{Au} has at most one solution, we consider $u\in H_0^2(\Omega)$ and obtain from the sesquilinear form that 
\[
	|u|_{H^2(\Omega)}\leq\frac{1}{\min(1-\mu, \mu)}A(u, u)=0, 
\]
which implies that $\frac{\partial^2 u}{\partial x^2_2}=0$. For any given $x_1\in [0, \Lambda]$, the following equation holds:
\[
\partial_{x_2} u(x_1, x_2)-\partial_{x_2} u(x_1, 0)=\int_0^{x_2} \partial^2_{tt} u(s, t){\rm d}t, 
\]
implying $\partial_{x_2} u(x_1, x_2)=\partial_{x_2} u(x_1, 0)=0$. By combining this result with the boundary condition on $S$, we conclude that
\[
u(x_1, x_2)-u(x_1, 0)=\int_0^{x_2} \partial_{t} u(s, t){\rm d}t=0,
\]
which shows $u(x_1, x_2)=u(x_1, 0)=0$. The lemma is proven through the application of the Fredholm alternative Theorem.
\end{proof}

The remainder of this subsection focuses on the stability estimate. Since the solution of \eqref{auxiliary} is quasi-periodic, it admits the Fourier series expansion
\begin{equation}\label{au-u}
	u(x_1, x_2)=\sum\limits_{n\in\mathbb{Z}} u_n(x_2) e^{{\rm i}\alpha_n x_1}.
\end{equation}
A simple calculation yields 
\begin{equation}\label{ODE}
	0=\Delta^2 u=\sum\limits_{n\in\mathbb{Z}}
	\left(u_n^{(4)}(x_2)-2\alpha_n^2 u_n^{(2)}(x_2)+\alpha_n^4 u_n(x_2)\right) e^{{\rm i}\alpha_n x_1},\quad x_2\in(0, h).
\end{equation}
If $\alpha_n\neq 0$, the general solution is given by 
\begin{equation}\label{un}
	u_n(x_2)=C_n^{(1)} e^{\alpha_n x_2}+C_n^{(2)} x_2 e^{\alpha_n x_2}+C_n^{(3)} e^{-\alpha_n x_2}+C_n^{(4)} x_2 e^{-\alpha_n x_2}.
\end{equation}
The Fourier coefficients $C_n^{(i)}$ satisfy the linear system
\begin{equation}\label{Coef}
\begin{bmatrix}
1 & 0 & 1 & 0\\
\alpha_n & 1 &-\alpha_n & 1\\
e^{\alpha_n h} & he^{\alpha_n h} & e^{-\alpha_n h} & he^{-\alpha_n h}\\
\alpha_n e^{\alpha_n h} & (1+\alpha_n h)e^{\alpha_n h} & -\alpha_n e^{-\alpha_n h} & (1-\alpha_n h)e^{-\alpha_n h}
\end{bmatrix}
\begin{bmatrix}
C_n^{(1)} \\ C_n^{(2)} \\ C_n^{(3)} \\C_n^{(4)}
\end{bmatrix}=
\begin{bmatrix}
0 \\ 0 \\ 0 \\ \varphi_n
\end{bmatrix}.
\end{equation}
It follows from a straightforward calculation that the solution of \eqref{Coef} is given by 
\begin{align*}
C_n^{(1)} &= \varphi_n h e^{\alpha_n h} \left(e^{2\alpha_n h}-1\right)/\Lambda_n,\\
C_n^{(2)} &= \varphi_n e^{\alpha_n h}\left(2\alpha_n h-e^{2\alpha_n h}+1\right)/\Lambda_n,\\
C_n^{(3)} &= -\varphi_n he^{\alpha_n h}\left(e^{2\alpha_n h}-1\right)/\Lambda_n,\\
C_n^{(4)} &=-\varphi_n e^{\alpha_n h}\left(2\alpha_n he^{2\alpha_n h}-e^{2\alpha_n h}+1\right)/\Lambda_n,
\end{align*}
where
\[
 \Lambda_n = 2e^{2\alpha_n h}-e^{4\alpha_n h}+4\alpha_n^2 h^2 e^{2\alpha_n h}-1. 
\]

In order to demonstrate that the solution is well-defined, it is necessary to verify that $\Lambda_n\neq 0$ for any $h>0$. Let
\[
f(t)=2e^{2t}-e^{4t}+4t^2 e^{2t}-1.
\]
It is easy to check that $ f(0)=0$ and
\[
	 f'(t)=4e^{2t}-4e^{4t}+8te^{2t}+8t^2e^{2t}=4e^{2t}\left(1-e^{2t}+2t+2t^2\right), \quad f'(0)=0,
\]
Let $g(t)=1-e^{2t}+2t+2t^2$. It follows that
\[
\left\{
\begin{aligned}
	&g'(t)=-2e^{2t}+2+4t, \quad g'(0)=0,\\
	&g''(t)=-4e^{2t}+4,
\end{aligned}
\right.
\]
Thus, it can be observed that:  
\begin{itemize}
\item[(i)] when $t>0$, $g''(t)<0$, $g'(t)$ decreases and is negative. Then $f'$ decreases and is negative, leading to a decrease and negative value in $f$. Thus $\Lambda_n<0$ for $\alpha_n h<0$.

\item[(ii)] when $t<0$, $g''(t)>0$, $g'(t)$ increases and is positive. Then $f'$ increases and is positive, 
resulting in an increase and positive value in $f$. Thus $\Lambda_n>0$ for $\alpha_n h>0$.
\end{itemize}

Combining the above discussion, $u_n$ is well-defined for arbitrary $h>0$. Furthermore, since $\Lambda_n\neq 0$,
we can also show that the auxiliary problem \eqref{auxiliary} is well-posed.

If $\alpha_n=0$, the governing equation \eqref{ODE} reduces to $\alpha_n^{(4)}(x_2)=0$. Hence, we obtain $u_n(x_2)=\left(h^{-2} x_2^3- h^{-1} x_2^2\right)\varphi_n,$ which shows that the coefficient $u_n$ is well-defined. 

Based on \eqref{au-u}, a straightforward computation yields 
\begin{align*}
&\partial_{x_1} u=\sum\limits_{n\in\mathbb{Z}} {\rm i}\alpha_n u_n(x_2) e^{{\rm i}\alpha_n x_1},\quad
\partial_{x_2} u=\sum\limits_{n\in\mathbb{Z}} u'_n(x_2) e^{{\rm i}\alpha_n x_1},\\
&\partial^2_{x_1} u=-\sum\limits_{n\in\mathbb{Z}} \alpha^2_n u_n(x_2) e^{{\rm i}\alpha_n x_1},\quad
\partial^2_{x_2} u=\sum\limits_{n\in\mathbb{Z}} u''_n(x_2) e^{{\rm i}\alpha_n x_1},\quad
\partial^2_{x_1 x_2} u=\sum\limits_{n\in\mathbb{Z}} {\rm i}\alpha_n u'_n(x_2) e^{{\rm i}\alpha_n x_1}.
\end{align*}
The $L^2$-norm is defined by
\begin{align*}
\|u\|^2_{L^2(\Omega)} = \int_0^h \int_0^\Lambda \left|\sum\limits_{n\in\mathbb{Z}} u_n(x_2) 
e^{{\rm i}\alpha_n x_1}\right|^2 {\rm d}x
=\int_0^h \Lambda\sum\limits_{n\in\mathbb{Z}} |u_n(x_2)|^2 {\rm d}x_2,
\end{align*}
the $H^1$-semi norm is
\begin{align*}
\left|u\right|^2_{H^1(\Omega)} &=\int_{\Omega} \left(\left|\frac{\partial u}{\partial x_1}\right|^2
	+\left|\frac{\partial u}{\partial x_2}\right|^2\right) {\rm d}x\\
&= \int_0^h \int_0^\Lambda \left|\sum\limits_{n\in\mathbb{Z}} {\rm i}\alpha_n u_n(x_2) e^{{\rm i}\alpha_n x_1}\right|^2
	+\left|\sum\limits_{n\in\mathbb{Z}} u'_n(x_2) e^{{\rm i}\alpha_n x_1}\right|^2 {\rm d}x \\
	&= \int_0^h \Lambda\sum\limits_{n\in\mathbb{Z}}\left(\alpha_n^2 |u_n(x_2)|^2+|u'_n(x_2)|^2\right) {\rm d}x_2,
\end{align*}
and the $H^2$-semi norm is
\begin{align}\label{H2-sn}
\left|u\right|^2_{H^2(\Omega)} &=
\int_{\Omega} \left(\left|\frac{\partial^2 u}{\partial x_1^2}\right|^2
	+\left|\frac{\partial^2 u}{\partial x_2^2}\right|^2+2\left|\frac{\partial^2 u}{\partial x_1 x_2}\right|^2\right) {\rm d}x\notag\\
&=  \int_0^h \int_0^\Lambda 
\left|-\sum\limits_{n\in\mathbb{Z}} \alpha^2_n u_n(x_2) e^{{\rm i}\alpha_n x_1}\right|^2
+\left|\sum\limits_{n\in\mathbb{Z}} u''_n(x_2) e^{{\rm i}\alpha_n x_1}\right|^2
+2\left|\sum\limits_{n\in\mathbb{Z}} {\rm i}\alpha_n u'_n(x_2) e^{{\rm i}\alpha_n x_1}\right|^2 {\rm d}x\notag\\
&= \int_0^h \Lambda\sum\limits_{n\in\mathbb{Z}}
\left(\alpha^4_n|u_n(x_2)|^2+|u_n''(x_2)|^2+2\alpha_n^2 |u'_n(x_2)|^2\right){\rm d}x_2.
\end{align}
Specially, when $\alpha_n=0$, it holds that
\[
	\int_0^h \left|u_n(x_2)\right|^2{\rm d}x_2<4h^4 \int_0^h \left|u_n''(x_2)\right|^2{\rm d}x_2.
\]

It is clear to note that the $H^2$-norm is bounded by the $H^2$-semi norm. Thus it suffices to show that 
\[
\left|u\right|^2_{H^2(\Omega)} \leq \|\varphi\|^2_{H^{1/2}(\Gamma)}.
\]
which requires estimating the three terms on the right-hand side of \eqref{H2-sn}.

First, we estimate the term involving $u_n$. It follows from a direct computation that 
\begin{align*}
C_n^{(1)}+C_n^{(2)}x_2 &= \frac{\varphi_n}{\Lambda_n}\left[(h-x_2)e^{3\alpha_n h}+(x_2-h+2\alpha_n h x_2)e^{\alpha_n h}\right],\\
C_n^{(3)}+C_n^{(4)}x_2 &=-\frac{\varphi_n}{\Lambda_n}\left[(h-x_2+2\alpha_n h x_2)e^{3\alpha_n h}+(x_2-h)e^{\alpha_n h}\right],
\end{align*}
and
\begin{align}
|u_n(x_2)|^2 &\leq 2\left|C_n^{(1)}+C_n^{(2)}x_2\right|^2 e^{2\alpha_n x_2}+ 2\left|C_n^{(3)}+C_n^{(4)}x_2\right|^2 e^{-2\alpha_n x_2}
\notag\\
&\leq\frac{4|\varphi_n|^2}{\Lambda_n^2}\Bigg[\left((h-x_2)^2 e^{6\alpha_n h}+(x_2-h+2\alpha_n h x_2)^2 e^{2\alpha_n h}\right)e^{2\alpha_n x_2}\notag\\
&\quad+\left((h-x_2+2\alpha_n h x_2)^2 e^{6\alpha_n h}+(x_2-h)^2 e^{2\alpha_n h}\right)e^{-2\alpha_n x_2}
\Bigg].\label{un}
\end{align}
For any given positive $h$, there exists a constant $C$ depending on $h$ such that
\[
\left\{
\begin{aligned}
	&\frac{1}{|\Lambda_n|}\leq C e^{-4\alpha_n h},\quad && \alpha_n>0\\
	&\frac{1}{|\Lambda_n|}\leq C, &&\alpha_n< 0
\end{aligned}
\right.
\]
Combining the above estimates gives 
\begin{eqnarray*}
|u_n(x_2)|^2\leq
\left\{
\begin{aligned}
	& C_1 |\varphi_n|^2\Big[(h-x_2)^2 e^{-2\alpha_n h+2\alpha_n x_2}\\
	& \qquad +(h-x_2 +2\alpha_n h x_2)^2 e^{-2\alpha_n h-2\alpha_n x_2}\Big],\, && \alpha_n>0,\\
	& C_1  |\varphi_n|^2\Big[(h-x_2)^2 e^{2\alpha_n h-2\alpha_n x_2}\\
	&\qquad +(x_2 -h+2\alpha_n h x_2)^2 e^{2\alpha_n h+2\alpha_n x_2}\Big],&&\alpha_n< 0,
\end{aligned}
\right.
\end{eqnarray*}

For $\alpha_n>0$, we deduce 
\begin{align*}
\int_0^h |u_n(x_2)|^2{\rm d}x_2 &\leq C_1 |\varphi_n|^2\Bigg[\frac{1}{4\alpha_n^3}
	+\frac{h^2\alpha_n+h}{\alpha_n^2}e^{-2\alpha_n h}\\
&\quad	-\left(2\alpha_n h^4+2h^3+\frac{4\alpha_n h-8h^2\alpha_n^2-1}{4\alpha_n^3}\right)e^{-4\alpha_n h}
	\Bigg]\\
	&\leq C_1 |\varphi_n|^2\left[\frac{1}{4\alpha_n^3}
	+\frac{h^2\alpha_n+h}{\alpha_n^2}e^{-2\alpha_n h}
	+\frac{8h^2\alpha_n^2+1}{4\alpha_n^3}e^{-4\alpha_n h}\right]\\
	&\leq C_2 \frac{ |\varphi_n|^2}{\alpha_n^3}
\end{align*}
When $\alpha_n<0$, it holds that
\begin{align*}
\int_0^h |u_n(x_2)|^2{\rm d}x_2 \leq C_2 \frac{|\varphi_n|^2}{|\alpha_n|^3},
\end{align*}
which gives 
\begin{align*}
\Lambda \sum\limits_{n\in\mathbb{Z}} \int_0^h
\alpha^4_n|u_n(x_2)|^2{\rm d}x_2 \leq C_2\Lambda
 \sum\limits_{n\in\mathbb{Z}} \alpha^4_n  |\varphi_n|^2\frac{1}{|\alpha_n|^3} \leq C_3\Lambda\|\varphi\|^2_{H^{1/2}(\Gamma)},
\end{align*}
where the positive constants $C_1, C_2$ and $C_3$ depend only on $h$.

Next is to estimate the term involving $u'_n$. A direct computation gives 
\begin{align*}
 u'_n(x_2)=\left[(C_n^{(1)}+C_n^{(2)} x_2)\alpha_n+C_n^{(2)}\right]e^{\alpha_n x_2}
+\left[-(C_n^{(3)}+C_n^{(4)} x_2)\alpha_n+C_n^{(4)}\right]e^{-\alpha_n x_2},
\end{align*}
which implies  
\begin{align*}
\left|u'_n(x_2)\right|^2 &\leq 4 \left(\left|C_n^{(1)}+C_n^{(2)} x_2\right|^2\alpha^2_n+\left|C_n^{(2)}\right|^2\right) e^{2\alpha_n x_2}
+4\left(\left|C_n^{(3)}+C_n^{(4)} x_2\right|^2\alpha^2_n+\left|C_n^{(4)}\right|^2\right)e^{-2\alpha_n x_2}\\
& \leq  4\alpha_n^2 \left[\left|C_n^{(1)}+C_n^{(2)} x_2\right|^2 e^{2\alpha_n x_2}
+\left|C_n^{(3)}+C_n^{(4)} x_2\right|^2 e^{-2\alpha_n x_2}\right]\\
&\quad 
+4\left(\left|C_n^{(2)}\right|^2 e^{2\alpha_n x_2}+\left|C_n^{(4)}\right|^2e^{-2\alpha_n x_2}\right).
\end{align*}

There exists a positive constant $C_4$ depending only on $h$ such that
\begin{align*}
|C_n^{(2)}|^2\leq C_4 e^{-2|\alpha_n| h}|\varphi_n|^2,\quad
|C_n^{(4)}|^2\leq C_4 e^{-2|\alpha_n| h}|\varphi_n|^2,
\end{align*}
and
\begin{align*}
\int_0^h \left(\left|C_n^{(2)}\right|^2 e^{2\alpha_n x_2}+\left|C_n^{(4)}\right|^2e^{-2\alpha_n x_2}\right) {\rm d}x_2
& \leq \frac{C_4}{2|\alpha_n|}\left(e^{2\alpha_n h}-1\right)e^{-2|\alpha_n| h}|\varphi_n|^2\\
&\quad +\frac{C_4}{2|\alpha_n|}\left(1-e^{-2\alpha_n h}\right)e^{-2|\alpha_n| h}|\varphi_n|^2,
\end{align*}
which imply
\begin{align*}
\Lambda \sum\limits_{n\in\mathbb{Z}}\int_0^h
\alpha^2_n|u'_n(x_2)|^2{\rm d}x_2 \leq \Lambda
\sum\limits_{n\in\mathbb{Z}} \alpha_n^2
\left( C_3\alpha_n^2\frac{|\varphi_n|^2}{|\alpha_n|^3}
+\frac{2C_4}{|\alpha_n|}  |\varphi_n|^2 \right)
\leq C_5\Lambda\|\varphi\|^2_{H^{1/2}(\Gamma)},
\end{align*}
where the positive constant $C_5$ depends only on $h$.

Finally, we estimate the term involving $u''_n$. A straightforward computation leads to 
\begin{align*}
u''_n(x_2)&=\left[C_n^{(1)}\alpha_n^2+2\alpha_n C_n^{(2)}+\alpha_n^2 x_2 C_n^{(2)}\right]e^{\alpha_n x_2}
+\left[C_n^{(3)}\alpha_n^2-2\alpha_n C_n^{(4)}+\alpha_n^2 x_2 C_n^{(4)}\right]e^{-\alpha_n x_2}\\
&=\alpha_n^2\left[\left(C_n^{(1)}+C_n^{(2)} x_2\right)e^{\alpha_n x_2}+\left(C_n^{(3)}+C_n^{(4)} x_2\right)e^{-\alpha_n x_2}\right] +2\alpha_n \left(C_n^{(2)} e^{\alpha_n x_2}-C_n^{(4)} e^{-\alpha_n x_2}\right)
\end{align*}
and
\begin{align*}
\left|u''_n(x_2)\right|^2 &\leq 4\alpha_n^4 \left[
\left|C_n^{(1)}+C_n^{(2)} x_2\right|^2 e^{2\alpha_n x_2}+\left|C_n^{(3)}+C_n^{(4)} x_2\right|^2 e^{-2\alpha_n x_2}\right]\\
&\quad +8\alpha_n^2\left[\left|C_n^{(2)}\right|^2 e^{2\alpha_n x_2}+\left|C_n^{(4)}\right|^2 e^{-2\alpha_n x_2}\right].
\end{align*}
When $\alpha_n\neq 0$, we have 
\begin{align*}
\int_0^h |u''_n(x_2)|^2{\rm d}x_2
\leq C 
\left( 4C_3\alpha_n^4\frac{  |\varphi_n|^2}{|\alpha_n|^3}
+16C_4\frac{\alpha_n^2}{|\alpha_n|}  |\varphi_n|^2 \right)
\leq C_6\left(1+\alpha_n^2\right)^{1/2} |\varphi_n|^2.
\end{align*}
When $\alpha_n=0$, it holds that 
\[
	\int_0^h \left|u_n''(x_2)\right|^2{\rm d}x_2=\frac{4}{h}\left|\varphi_n\right|^2=\frac{4}{h}\left(1+\alpha_n^2\right)^{1/2}\left|\varphi_n\right|^2,
\]
Hence, we deduce that 
\[
	\Lambda \sum\limits_{n\in\mathbb{Z}}\int_0^h |u''_n(x_2)|^2{\rm d}x_2
	\leq C_7\sum\limits_{n\in\mathbb{Z}} \left(1+\alpha_n^2\right)^{1/2}\left|\varphi_n\right|^2
	=C_7\|\varphi\|^2_{H^{1/2}(\Gamma)},
\]
where the positive constant $C_6$ depends only on $h$.

Combining Lemma \ref{lemmaau} with the preceding discussion, we demonstrate that the
 auxiliary problem \eqref{auxiliary} is uniquely solvable. Moreover, the solution satisfies the estimate
\[
	\|u\|_{H^2(\Omega)}\leq C\Lambda\|\varphi\|_{H^{1/2}(\Gamma)},
\]
where the positive constant $C$ depends only on $h$.

\end{document}